\documentclass[11pt]{amsart}
\usepackage{mathrsfs}
\usepackage{amssymb,latexsym}
\usepackage{pstcol,pstricks,color}

 \numberwithin{equation}{section}

\newtheorem{theorem}{Theorem}[section]

\newtheorem{corollary}[theorem]{Corollary}

\newtheorem{remark}[theorem]{Remark}

\newtheorem{example}[theorem]{Example}

\def\qed{\hfill $\Box$}
\def\pf{\noindent {\it Proof.} }

\title{ Weighted Dyck paths with special restrictions on the levels of valleys  }

\begin{document}
\maketitle
\begin{center}
Yidong Sun$^{\dag}$\footnote{Corresponding author: Yidong Sun.}, Qianqian Liu$^{\ddag}$ and Yanxin Liu$^{\sharp}$

School of Science, Dalian Maritime University, 116026 Dalian, P.R. China\\[5pt]

{\it Emails: $^{\dag}$sydmath@dlmu.edu.cn, $^{\ddag}$lqq961106@dlmu.edu.cn, $^{\sharp}$lyx@dlmu.edu.cn }

\end{center}\vskip0.2cm

\subsection*{Abstract} This paper concentrates on the set $\mathcal{V}_n$ of weighted Dyck paths of length $2n$ with special restrictions on the level of valleys. We first give its explicit formula of the counting generating function in terms of certain weight functions. When the weight functions are specialized, some connections are builded between $\mathcal{V}_n$ and other classical combinatorial structures such as $(a,b)$-Motzkin paths, $q$-Schr\"{o}der paths, Delannoy paths and complete $k$-ary trees. Some bijections are also established between these settings and $\mathcal{V}_n$ subject to certain special weight functions.

\medskip

{\bf Keywords}: Dyck path, $(a, b)$-Motzkin path, Schr\"{o}der path, Delannoy path, Narayana polynomial, Fuss-Catalan number.

\noindent {\sc 2020 Mathematics Subject Classification}: Primary 05A15; Secondary 05A10, 05A19.

{\bf \section{ Introduction } }

A {\it Dyck path} of length $2n$ is a lattice path from $(0, 0)$ to $(2n, 0)$ in the first quadrant of the xy-plane that consists of up steps $\mathbf{u}=(1, 1)$ and down steps $\mathbf{d}=(1, -1)$. See \cite[p.204]{StanleyEC} and \cite{Deutsch99}.

A point of a Dyck path with ordinate $\ell$ is said to be at {\it level} $\ell$. A step of a Dyck path is said to be at level $\ell$ if the ordinate of its endpoint is $\ell$. A {\it peak (valley)} in a Dyck path is an occurrence of $\mathbf{ud}$ ($\mathbf{du}$). By the {\it level of a peak (valley)} we mean the level of the intersection point of its two steps.
A {\it pyramid} in a Dyck path is a section of the form $\mathbf{u}^{h}\mathbf{d}^{h}$, a succession of $h$ up steps followed immediately by $h$ down steps, where $h$ is called the {\it height} of the pyramid. A pyramid $\mathbf{u}^{h}\mathbf{d}^{h}$ is {\it maximal} if it cannot be extended to a pyramid $\mathbf{u}^{h+1}\mathbf{d}^{h+1}$. The {\it altitude} of a maximal pyramid is the level of its last $\mathbf{d}$-step.

By a {\it return step} we mean a $\mathbf{d}$-step at level $0$. Dyck paths that have exactly one return step are said to be {\it primitive}. If $\mathbf{P}_1$ and $\mathbf{P}_2$ are Dyck paths, then we define $\mathbf{P}_1\mathbf{P}_2$ as the concatenation of $\mathbf{P}_1$ and $\mathbf{P}_2$, and define $\hat{\mathbf{P}}_1=\mathbf{u}\mathbf{P}_1\mathbf{d}$ as the elevation of $P_1$. Naturally, $\hat{\mathbf{P}}_1$ is primitive.

Let $\varepsilon$ be the empty path, that is a dot path. If $\mathcal{P}_1$ and $\mathcal{P}_2$ are sets of Dyck paths, we define the {\it concatenation} $\mathcal{P}_1\mathcal{P}_2$ of $\mathcal{P}_1$ and $\mathcal{P}_2$ by
$$\mathcal{P}_1\mathcal{P}_2=\{\mathbf{P}_1\mathbf{P}_2| \mathbf{P}_1\in \mathcal{P}_1, \mathbf{P}_2\in \mathcal{P}_2\} $$
and the {\it elevation} $\hat{\mathcal{P}}_1$ of $\mathcal{P}_1$ by
$$\hat{\mathcal{P}}_1=\{ \hat{\mathbf{P}}_1| \mathbf{P}_1\in \mathcal{P}_1\}. $$
Clearly, $\mathcal{P}_1\{\varepsilon\}=\{\varepsilon\}\mathcal{P}_1=\mathcal{P}_1.$

Let $\mathcal{D}_n$ denote the set of all Dyck paths of length $2n$. It is well known \cite{Stanley} that $\mathcal{D}_n$ is counted by Catalan numbers $C_n=\frac{1}{n+1}\binom{2n}{n}$, which has generating function $C(x)=\frac{1-\sqrt{1-4x}}{2x}$ satisfying the relation
$$C(x)=1+xC(x)^2=\frac{1}{1-xC(x)}. $$

Barcucci et at. \cite{Barcucci} introduced the concept, {\it non-decreasing Dyck path}, a Dyck path where the levels of its valleys form a non-decreasing sequence along the path, and since then several papers dedicated to the topic have appeared \cite{Czabarka1, Czabarka2, Florez1}, some bijections have been found between non-decreasing Dyck paths and other classical combinatorial settings such as directed column-convex polyominoes, Elena trees \cite{DeutschProd, Prodinger}. The same concept has been extended to Motzkin paths \cite{Florez2} and to t-Dyck paths \cite{Florez3}. Motivated by the work of Barcucci et at. \cite{Barcucci}, Czabarka et at. \cite{Czabarka1, Czabarka2} and Fl\'{o}rez et at. \cite{Florez1}, in this paper we consider a class of weighted Dyck paths with another restriction on the levels of valleys.

Let $\mathcal{A}_{n}$ denote the set of all nonzero-weighted primitive Dyck paths $\mathbf{P}$ of length $2n$ such that the valleys (if exist) of $\mathbf{P}$ have the same level. The weight of each part of $\mathbf{P}$ is assigned as follows:

$1)$ A maximal pyramid of height $k$ at altitude $0$ is weighted by $\gamma_k$ for $k\geq 1$;

$2)$ A maximal pyramid of height $k$ at altitude $\geq 1$ is weighted by $\alpha_k$ for $k\geq 1$;

$3)$ The beginning segment $\mathbf{u}^k$ of $\mathbf{P}$ is weighted by $\beta_k$ if its valleys are at level $k$ for $k\geq 1$.

The {\it weight} of $\mathbf{P}$, denoted by $w(\mathbf{P})$, is the product of the weight of each part of $\mathbf{P}$. For example, $w(\mathbf{u}^k\mathbf{d}^k)=\gamma_k$ and $w(\mathbf{u}^k\mathbf{u}^i\mathbf{d}^i\mathbf{u}^j\mathbf{d}^j\mathbf{d}^k)=\beta_k\alpha_i\alpha_j$ for $k, i, j\geq 1$. The {\em weight} of $\mathcal{A}_{n}$, denoted by $w(\mathcal{A}_{n})$, is the sum of the total weights of all paths in $\mathcal{A}_{n}$. Note that the weight of each $\mathbf{P}\in \mathcal{A}_{n}$ is nonzero.

Let $\mathcal{A}=\bigcup_{n=1}^{\infty}\mathcal{A}_{n}$, define
$$\mathcal{V}=\{\varepsilon\}+\mathcal{A}+\mathcal{A}\mathcal{A}+\mathcal{A}\mathcal{A}\mathcal{A}+\cdots. $$

Let $\mathcal{V}_n$ denote the set of weighted Dyck paths $\mathbf{P}$ of length $2n$ in $\mathcal{V}$. Then the valleys (if exists) in each primitive part of $\mathbf{P}$ have the same level. See Figure 1 as an example.

\begin{figure}[h] \setlength{\unitlength}{0.5mm}

\begin{center}
\begin{pspicture}(14.5,4)
\psset{xunit=15pt,yunit=15pt}\psgrid[subgriddiv=1,griddots=4,
gridlabels=4pt](0,0)(28,7)

\psline(0,0)(6,6)(9,3)(10,4)(11,3)(12,4)(16,0)(19,3)(20,2)(21,3)(24,0)(26,2)(28,0)

\pscircle*(0,0){0.06}\pscircle*(1,1){0.06}\pscircle*(2,2){0.06}
\pscircle*(3,3){0.06}\pscircle*(4,4){0.06}\pscircle*(5,5){0.06}
\pscircle*(6,6){0.06}\pscircle*(7,5){0.06}\pscircle*(8,4){0.06}
\pscircle*(9,3){0.06}\pscircle*(10,4){0.06}\pscircle*(11,3){0.06}
\pscircle*(12,4){0.06}\pscircle*(13,3){0.06}\pscircle*(14,2){0.06}
\pscircle*(15,1){0.06}\pscircle*(16,0){0.06}\pscircle*(17,1){0.06}
\pscircle*(18,2){0.06}\pscircle*(19,3){0.06}

\pscircle*(28,0){0.06}\pscircle*(27,1){0.06}\pscircle*(26,2){0.06}
\pscircle*(25,1){0.06}\pscircle*(24,0){0.06}\pscircle*(23,1){0.06}
\pscircle*(22,2){0.06}\pscircle*(21,3){0.06}\pscircle*(20,2){0.06}

\end{pspicture}
\end{center}

\caption{\small An example $\mathbf{P}=\mathbf{u}^3(\mathbf{u}^3\mathbf{d}^3)(\mathbf{u}\mathbf{d})^2\mathbf{d}^3\mathbf{u}^2(\mathbf{u}\mathbf{d})^2\mathbf{d}^2\mathbf{u}^2\mathbf{d}^2\in \mathcal{V}_{14}$ with the weight $w(\mathbf{P})=\alpha_1^4\alpha_3\beta_2\beta_3\gamma_2$. }

\end{figure}

Set $V_n=w(\mathcal{V}_n)$ and define the generating function $V_{\alpha, \beta, \gamma}(x)$ for $V_n$ related to the weight functions $\alpha(x)=\sum_{k\geq 1}\alpha_kx^{k}$, $\beta(x)=\sum_{k\geq 1}\beta_kx^{k}$ and $\gamma(x)=\sum_{k\geq 1}\gamma_kx^{k}$, i.e.,
$$V_{\alpha, \beta, \gamma}(x)=\sum_{n=0}^{\infty}V_nx^n. $$

In the paper, we concentrate on the generating function $V_{\alpha, \beta, \gamma}(x)$ and build some connections between $\mathcal{V}_n$ and other classical combinatorial settings such as $(a,b)$-Motzkin paths, $q$-Schr\"{o}der paths, Delannoy paths and complete $k$-ary trees. Precisely, the next section gives the explicit formula for $V_{\alpha, \beta, \gamma}(x)$ in terms of the weight functions $\alpha(x)$, $\beta(x)$ and $\gamma(x)$. In the third section, several special examples are considered when the weight functions $\alpha(x), \beta(x)$ and $\gamma(x)$ are specialized. Some bijections are also established between these settings and $\mathcal{V}_n$ subject to certain special weight functions $\alpha(x), \beta(x)$ and $\gamma(x)$.

\vskip0.5cm
\section{The generating function for $V_n$}

In this section, we derive the generating function $V_{\alpha, \beta, \gamma}(x)$ for $V_n$ or for $\mathcal{V}$ in terms of the weight functions $\alpha(x)$ , $\beta(x)$ and $\gamma(x)$.

\begin{theorem}\label{theo2.1}
The generating function $V_{\alpha, \beta, \gamma}(x)$ for $V_n$ reads
\begin{eqnarray}\label{eqn 2.1}
V_{\alpha, \beta, \gamma}(x)=\sum_{n=0}^{\infty}V_nx^n=\frac{1}{1-\gamma(x)-\frac{\alpha(x)^2\beta(x)}{1-\alpha(x)}} .
\end{eqnarray}
\end{theorem}
\pf It is clear that $\mathcal{V}$ can be partitioned as follows
$$\mathcal{V}=\{\varepsilon\}+\mathcal{A}\mathcal{V}. $$

The empty path $\varepsilon$ contributes the weight 1. For any path $\mathbf{P}\in \mathcal{A}$, when $\mathbf{P}$ is primitive, there are two cases to be considered:

$1)$ If $\mathbf{P}$ is a maximal pyramid of height $h$ at altitude $0$, which has the weight $\gamma_h$ for $h\geq 1$, this case produces the weight function $\gamma(x)=\sum_{h\geq 1}\gamma_hx^{h}$;

$2)$ If $\mathbf{P}$ has $r\geq 2$ peaks such that the $r-1$ valleys have the level $k\geq 1$, so $\mathbf{P}$ is of form $\mathbf{P}=\mathbf{u}^k\mathbf{u}^{i_1}\mathbf{d}^{i_1}\cdots \mathbf{u}^{i_r}\mathbf{d}^{i_r}\mathbf{d}^{k}$, where $\mathbf{u}^{i_j}\mathbf{d}^{i_j}$ are maximal pyramids at altitude $k\geq 1$ and of hight $i_j\geq 1$ for $1\leq j\leq r$. Each kind of maximal pyramids gives the weight function $\alpha(x)=\sum_{i_j\geq 1}\alpha_{i_j}x^{i_j}$. Then these lead to the weight function $\beta_kx^k\alpha(x)^{r}$. Summarizing $k\geq 1$ and $r\geq 2$, we have the counting generating function $\frac{\alpha(x)^2\beta(x)}{1-\alpha(x)}$.

These two cases give that the weight function for $\mathcal{A}$ is
$$\gamma(x)+\frac{\alpha(x)^2\beta(x)}{1-\alpha(x)}. $$

Hence, we have the relation of the counting generating function for $\mathcal{V}$ as follows,
$$V_{\alpha, \beta, \gamma}(x)=1+\Big(\gamma(x)+\frac{\alpha(x)^2\beta(x)}{1-\alpha(x)}\Big)V_{\alpha, \beta, \gamma}(x). $$
When solve it for $V_{\alpha, \beta, \gamma}(x)$, one obtains (\ref{eqn 2.2}), the desired result.   \qed

\vskip0.2cm

When $\gamma(x)=\alpha(x)\beta(x)$, we write $V_{\alpha, \beta}(x)$ instead of $V_{\alpha, \beta, \alpha\beta}(x)$ for short. In this case, we have

\begin{corollary}\label{coro2.2}
\begin{eqnarray}\label{eqn 2.2}
V_{\alpha, \beta}(x)=\frac{1-\alpha(x)}{1-\alpha(x)-\alpha(x)\beta(x)} .
\end{eqnarray}
\end{corollary}

For example, when $\alpha(x)=\frac{x}{1-x}, \beta(x)=\frac{x}{1-2x}$, (\ref{eqn 2.2}) implies that
\begin{eqnarray*}
V_{\alpha, \beta}(x)=\frac{1-\frac{x}{1-x}}{1-\frac{x}{1-x}-\frac{x^2}{(1-x)(1-2x)}}=1+\frac{x^2}{(1-x)(1-3x)}=1+\sum_{n=1}^{\infty}\frac{1}{2}(3^{n-1}-1)x^{n} .
\end{eqnarray*}
Another example of (\ref{eqn 2.2}), when $\alpha(x)=\beta(x)=\frac{x}{1-x}$ leads to
\begin{eqnarray*}
V_{\alpha, \beta}(x)=\frac{1-\frac{x}{1-x}}{1-\frac{x}{1-x}-\frac{x^2}{(1-x)^2}}=1+\frac{x^2}{1-3x+x^2}=1+\sum_{n=1}^{\infty}F_{2(n-1)}x^{n} ,
\end{eqnarray*}
where $F_{k}$ is the Fibonacci sequence defined by $F_k=F_{k-1}+F_{k-2}$ with $F_0=0, F_1=1.$

\vskip0.5cm
\section{Several special cases for $\gamma(x)=\alpha(x)\beta(x)$}

In this section, we mainly consider $V_{\alpha, \beta, \gamma}(x)$ under the condition $\gamma(x)=\alpha(x)\beta(x)$. When the weight functions $\alpha(x), \beta(x)$ are specialized, several examples exhibit that $V_n$ are closely related to some classical sequences such as Catalan numbers, $(a, b)$-Motzkin numbers, $q$-Schr\"{o}der numbers, Narayana polynomials, central Delannoy numbers, Fibonacci numbers, Fuss-Catalan numbers and so on.

\subsection{The special case related to $(a, b)$-Motzkin paths}

Let
$$M^{(a, b)}(x)=\sum_{n\geq 0}M_n^{(a, b)}x^n=\frac{1-ax-\sqrt{(1-ax)^2-4bx^2}}{2bx^2}$$
be the generating function for $(a, b)$-Motzkin numbers $M^{(a, b)}_n$, which count the set $\mathcal{M}_n^{(a, b)}$ of $(a, b)$-Motzkin paths of length $n$ from $(0, 0)$ to $(n, 0)$ in the first quadrant of the xy-plane that consists of up steps $\mathbf{u}=(1, 1)$ with weight $1$, down steps $\mathbf{d}=(1, -1)$ with weight $b$ and horizontal steps $\mathbf{h}=(1, 0)$ with weight $a$ \cite{Woan}. When $(a, b)=(1, 1)$, we obtain Motzkin paths. Note that $M^{(a, b)}(x)$ satisfies the relations
$$M^{(a, b)}(x)=1+axM^{(a, b)}(x)+bx^2M^{(a, b)}(x)^2=\frac{1}{1-ax-bx^2M^{(a, b)}(x)}.$$

\begin{example}  The special case in Corollary \ref{coro2.2} when
$$\alpha(x)=ax, \beta(x)=\frac{b}{a}xM^{(a, b)}(x), \gamma(x)=bx^2M^{(a, b)}(x) $$
generates
\begin{eqnarray*}
V_{\alpha, \beta}(x)=\frac{1-ax}{1-ax-bx^2M^{(a, b)}(x)}=(1-ax)M^{(a, b)}(x),
\end{eqnarray*}
or equivalently, $V_n=M^{(a, b)}_n-aM^{(a, b)}_{n-1}$ for $n\geq 0$.
\end{example}

Let $\mathcal{X}_n$ denote the set of $(a, b)$-Motzkin paths of length $n$ such that the first step is not an $\mathbf{h}$-step. Then it is clear that $\mathcal{X}_n$ is counted by $M^{(a, b)}_n-aM^{(a, b)}_{n-1}$. This indicates the following corollary.
\begin{corollary}
There exists a simple bijection between the set $\mathcal{X}_n$ and the set $\mathcal{V}_n$ with the weight functions $\alpha(x)=ax, \beta(x)=\frac{b}{a}xM^{(a, b)}(x)$ and $\gamma(x)=bx^2M^{(a, b)}(x)$.
\end{corollary}
\pf It is trivial for $n=0$, an empty path, and for $n=1$, an empty set. According to the weight functions $\alpha(x)=ax, \beta(x)=\frac{b}{a}xM^{(a, b)}(x)$ and $\gamma(x)=bx^2M^{(a, b)}(x)$, for any $\mathbf{P}\in \mathcal{V}_n$ with $n\geq 2$, each primitive part of $\mathbf{P}$ has one of the following two forms: 1) it is a maximal pyramid of height $k\geq 2$ at altitude $0$ with weight $\gamma_{k-2}=bM_{k-2}^{(a,b)}$; 2) it is of form $\mathbf{u}^{k}(\mathbf{ud})^r\mathbf{d}^k$ with weight $\alpha_1^{r}\beta_{k-1}=a^r\cdot \frac{b}{a}M_{k-1}^{(a,b)}$ for $k\geq 1$ and $r\geq 2$ because the weight functions $\alpha(x)=ax$ requires that the height of all maximal pyramids at altitude $k\geq 1$ is equal to $1$. Then $\mathbf{P}$ can be written uniquely as $\mathbf{P}=\mathbf{P}_1\mathbf{P}_2$, where $\mathbf{P}_1$ is the first primitive part of length at least $4$. That is,  $\mathbf{P}_1=\mathbf{u}^{k}(\mathbf{ud})^r\mathbf{d}^k$ for certain $k, r\geq 1$.

When $r=1$, $\mathbf{P}_1=\mathbf{u}^{k+1}\mathbf{d}^{k+1}$ is a maximal pyramid at altitude 0 with weight $\gamma_{k-1}=bM^{(a, b)}_{k-1}$. Equivalently, we associate a primitive $(a, b)$-Motzkin path $\mathbf{u}\mathbf{Q}_{k-1}\mathbf{d}$ to $\mathbf{P}_1$, where $\mathbf{Q}_{k-1}\in \mathcal{M}_{k-1}^{(a, b)}$. In other words, we assign any $(a, b)$-Motzkin path $\mathbf{u}\mathbf{Q}_{k-1}\mathbf{d}$ to a maximal pyramid of height $k+1$ at altitude $0$.

When $r\geq 2$, $\mathbf{P}_1$ has weight $\alpha_1^{r}\beta_{k-1}=a^{r-1}bM^{(a, b)}_{k-1}$. Equivalently, a primitive $(a, b)$-Motzkin path $\mathbf{u}\mathbf{Q}_{k-1}\mathbf{d}$ is allocated to the first maximal $\mathbf{u}$-segment $\mathbf{u}^{k+1}$ in $\mathbf{P}_1$, and an $\mathbf{h}$ step with weight $a$ is assigned to each maximal pyramid of height 1 in $\mathbf{P}_1$ except for the first one. In other words,
we associate an $(a, b)$-Motzkin path $\mathbf{u}\mathbf{Q}_{k-1}\mathbf{d}\mathbf{h}^{r-1}$ to $\mathbf{P}_1$, where $\mathbf{Q}_{k-1}\in \mathcal{M}_{k-1}^{(a, b)}$.

Now we can recursively establish a bijection $\phi$ between $\mathcal{V}_n$ and $\mathcal{X}_n$ as follows. For any $\mathbf{P}=\mathbf{P}_1\mathbf{P}_2\in \mathcal{V}_n$, where $\mathbf{P}_1=\mathbf{u}^{k}(\mathbf{ud})^r\mathbf{d}^k$ for certain $k, r\geq 1$, if $\mathbf{P}_1$ is associated by an $(a, b)$-Motzkin path $\mathbf{u}\mathbf{Q}_{k-1}\mathbf{d}\mathbf{h}^{r-1}$ as above, then $\phi(\mathbf{P})=\mathbf{u}\mathbf{Q}_{k-1}\mathbf{d}\mathbf{h}^{r-1}\phi(\mathbf{P}_2)\in \mathcal{X}_n$.

Conversely, for any $(a, b)$-Motzkin path $\mathbf{M}\in\mathcal{X}_n$, $\mathbf{M}$ can be written uniquely as $\mathbf{M}=\mathbf{u}\mathbf{Q}_{k-1}\mathbf{d}\mathbf{h}^{r-1}\mathbf{M}_2$ for certain $k, r\geq 1$, $\mathbf{Q}_{k-1}\in \mathcal{M}_{k-1}^{(a, b)}$ and $\mathbf{M}_2\in \mathcal{X}_{n-k-r}$, the inverse $\phi^{-1}$ is built by $\phi^{-1}(\mathbf{M})=\mathbf{P}_1\phi^{-1}(\mathbf{M}_2)\in \mathcal{V}_n$, where $\mathbf{P}_1=\mathbf{u}^{k}(\mathbf{ud})^r\mathbf{d}^k$ is associated with an $(a, b)$-Motzkin path $\mathbf{u}\mathbf{Q}_{k-1}\mathbf{d}\mathbf{h}^{r-1}$. One can see that the total weight of all $\mathbf{P}_1$ is $a^{r-1}bM^{(a, b)}_{k-1}$.    \qed  \vskip0.2cm

In order to give a more intuitive view on the bijection $\phi$, we present a pictorial description of
$\phi$ for the case $\mathbf{P}=\mathbf{u}^5\mathbf{d}^5\mathbf{u}^3(\mathbf{u}\mathbf{d})^4\mathbf{d}^3\mathbf{u}^2\mathbf{d}^2$ with weight $bM^{(a, b)}_{3}\cdot\frac{b}{a}M^{(a, b)}_{2}a^{4}\cdot bM^{(a, b)}_{0}=a^3b^3(a^2+b)(a^3+3ab)$. Clearly, $\phi(\mathbf{P})=\mathbf{u}\mathcal{M}_{3}^{(a, b)}\mathbf{d}\mathbf{u}\mathcal{M}_{2}^{(a, b)}\mathbf{d}\mathbf{hhh}\mathbf{u}\mathbf{d}$ with $\mathbf{h}$ steps weighted by $a$ and $\mathbf{d}$ steps weighted by $b$. Precisely, if the first primitive part $\mathbf{u}^5\mathbf{d}^5$ of $\mathbf{P}$ is associated with a primitive $(a, b)$-Motzkin path $\mathbf{u}\mathbf{Q}_{3}\mathbf{d}$, the beginning $\mathbf{u}$-segment $\mathbf{u}^4$ of the second primitive part $\mathbf{u}^3(\mathbf{u}\mathbf{d})^4\mathbf{d}^3$ of $\mathbf{P}$ is associated with a primitive $(a, b)$-Motzkin path $\mathbf{u}\mathbf{Q}_{2}\mathbf{d}$ and the third primitive part $\mathbf{u}^2\mathbf{d}^2$ of $\mathbf{P}$ is associated with a primitive $(a, b)$-Motzkin path $\mathbf{u}\mathbf{Q}_{0}\mathbf{d}$, where $\mathbf{Q}_i\in \mathcal{M}_{i}^{(a, b)}$ for $i=2, 3$ and $\mathbf{Q}_0$ is an empty path, then $\phi(\mathbf{P})=\mathbf{u}\mathbf{Q}_3\mathbf{d}\mathbf{u}\mathbf{Q}_2\mathbf{d}\mathbf{hhh}\mathbf{u}\mathbf{d}\in \mathcal{X}_{14}$. See Figure 2.

\begin{figure}[h] \setlength{\unitlength}{0.5mm}

\begin{center}
\begin{pspicture}(15,3.5)
\psset{xunit=15pt,yunit=15pt}\psgrid[subgriddiv=1,griddots=4,
gridlabels=4pt](0,0)(29,7)

\psline(0,0)(5,5)(10,0)(14,4)(15, 3)(16,4)(17,3)(18,4)(19,3)(20,4)(24,0)(26,2)(28,0)

\pscircle*(0,0){0.06}\pscircle*(1,1){0.06}\pscircle*(2,2){0.06}
\pscircle*(3,3){0.06}\pscircle*(4,4){0.06}\pscircle*(5,5){0.06}
\pscircle*(6,4){0.06}\pscircle*(7,3){0.06}\pscircle*(8,2){0.06}
\pscircle*(9,1){0.06}\pscircle*(10,0){0.06}\pscircle*(11,1){0.06}
\pscircle*(12,2){0.06}\pscircle*(13,3){0.06}\pscircle*(14,4){0.06}
\pscircle*(15,3){0.06}\pscircle*(16,4){0.06}\pscircle*(17,3){0.06}
\pscircle*(18,4){0.06}\pscircle*(19,3){0.06}\pscircle*(20,4){0.06}

\pscircle*(28,0){0.06}\pscircle*(27,1){0.06}\pscircle*(26,2){0.06}
\pscircle*(25,1){0.06}\pscircle*(24,0){0.06}\pscircle*(23,1){0.06}
\pscircle*(22,2){0.06}\pscircle*(21,3){0.06}

\put(0.4,2.8){$w=bM^{(a, b)}_{3}=b(a^3+3ab)$} \put(6.5,2.5){$w=\frac{b}{a}M^{(a, b)}_{2}a^{4}=a^3b(a^2+b)$} \put(12.5,1.4){$w=bM^{(a, b)}_{0}=b$}
\put(6.5,0.1){$\mathbf{P}=\mathbf{u}^5\mathbf{d}^5\mathbf{u}^3(\mathbf{u}\mathbf{d})^4\mathbf{d}^3\mathbf{u}^2\mathbf{d}^2$}

\end{pspicture}
\end{center}\vskip0.5cm

$\Updownarrow \phi$
\vskip1.2cm

\begin{center}
\begin{pspicture}(15,7.5)
\psset{xunit=15pt,yunit=15pt}\psgrid[subgriddiv=1,griddots=4,
gridlabels=4pt](0,0)(29,16)

\psline(0,0)(1,1)(4,1)(5,0)(6,1)(8,1)(9,0)(12, 0)(13,1)(14,0)
\psline(0,3)(2,5)(3,4)(4,4)(5,3)(6,4)(8,4)(9,3)(12, 3)(13,4)(14,3)
\psline(0,6)(1,7)(2,7)(3,8)(5,6)(6,7)(8,7)(9,6)(12, 6)(13,7)(14,6)
\psline(0,9)(2,11)(3,11)(4,10)(5,9)(6,10)(8,10)(9,9)(12, 9)(13,10)(14,9)

\psline(15,0)(16,1)(19,1)(20,0)(22,2)(23,1)(24,0)(27,0)(28,1)(29,0)
\psline(15,3)(17,5)(18,4)(19,4)(20,3)(22,5)(23,4)(24,3)(27,3)(28,4)(29,3)
\psline(15,6)(16,7)(17,7)(18,8)(20,6)(22,8)(23,7)(24,6)(27,6)(28,7)(29,6)
\psline(15,9)(17,11)(18,11)(19,10)(20,9)(22,11)(23,10)(24,9)(27,9)(28,10)(29,9)

\psline(7,13)(8,14)(11,14)(12,13)(13,14)(15,14)(16,13)(19,13)(20,14)(21,13)
\pscircle*(7,13){0.06}\pscircle*(8,14){0.06}\pscircle*(11,14){0.06}
\pscircle*(12,13){0.06}\pscircle*(13,14){0.06}\pscircle*(15,14){0.06}
\pscircle*(16,13){0.06}\pscircle*(17,13){0.06}\pscircle*(18,13){0.06}
\pscircle*(19,13){0.06}\pscircle*(20,14){0.06}\pscircle*(21,13){0.06}

\pscircle*(0,0){0.06}\pscircle*(1,1){0.06}\pscircle*(2,1){0.06}
\pscircle*(3,1){0.06}\pscircle*(4,1){0.06}\pscircle*(5,0){0.06}
\pscircle*(6,1){0.06}\pscircle*(7,1){0.06}\pscircle*(8,1){0.06}
\pscircle*(9,0){0.06}\pscircle*(10,0){0.06}\pscircle*(11,0){0.06}
\pscircle*(12,0){0.06}\pscircle*(13,1){0.06}\pscircle*(14,0){0.06}

\pscircle*(0,3){0.06}\pscircle*(1,4){0.06}\pscircle*(2,5){0.06}
\pscircle*(3,4){0.06}\pscircle*(4,4){0.06}\pscircle*(5,3){0.06}
\pscircle*(6,4){0.06}\pscircle*(7,4){0.06}\pscircle*(8,4){0.06}
\pscircle*(9,3){0.06}\pscircle*(10,3){0.06}\pscircle*(11,3){0.06}
\pscircle*(12,3){0.06}\pscircle*(13,4){0.06}\pscircle*(14,3){0.06}

\pscircle*(0,6){0.06}\pscircle*(1,7){0.06}\pscircle*(2,7){0.06}
\pscircle*(3,8){0.06}\pscircle*(4,7){0.06}\pscircle*(5,6){0.06}
\pscircle*(6,7){0.06}\pscircle*(7,7){0.06}\pscircle*(8,7){0.06}
\pscircle*(9,6){0.06}\pscircle*(10,6){0.06}\pscircle*(11,6){0.06}
\pscircle*(12,6){0.06}\pscircle*(13,7){0.06}\pscircle*(14,6){0.06}

\pscircle*(0,9){0.06}\pscircle*(1,10){0.06}\pscircle*(2,11){0.06}
\pscircle*(3,11){0.06}\pscircle*(4,10){0.06}\pscircle*(5,9){0.06}
\pscircle*(6,10){0.06}\pscircle*(7,10){0.06}\pscircle*(8,10){0.06}
\pscircle*(9,9){0.06}\pscircle*(10,9){0.06}\pscircle*(11,9){0.06}
\pscircle*(12,9){0.06}\pscircle*(13,10){0.06}\pscircle*(14,9){0.06}

\pscircle*(15,0){0.06}\pscircle*(16,1){0.06}\pscircle*(17,1){0.06}
\pscircle*(18,1){0.06}\pscircle*(19,1){0.06}\pscircle*(20,0){0.06}
\pscircle*(21,1){0.06}\pscircle*(22,2){0.06}\pscircle*(23,1){0.06}
\pscircle*(24,0){0.06}\pscircle*(25,0){0.06}\pscircle*(26,0){0.06}
\pscircle*(27,0){0.06}\pscircle*(28,1){0.06}\pscircle*(29,0){0.06}

\pscircle*(15,3){0.06}\pscircle*(16,4){0.06}\pscircle*(17,5){0.06}
\pscircle*(18,4){0.06}\pscircle*(19,4){0.06}\pscircle*(20,3){0.06}
\pscircle*(21,4){0.06}\pscircle*(22,5){0.06}\pscircle*(23,4){0.06}
\pscircle*(24,3){0.06}\pscircle*(25,3){0.06}\pscircle*(26,3){0.06}
\pscircle*(27,3){0.06}\pscircle*(28,4){0.06}\pscircle*(29,3){0.06}

\pscircle*(15,6){0.06}\pscircle*(16,7){0.06}\pscircle*(17,7){0.06}
\pscircle*(18,8){0.06}\pscircle*(19,7){0.06}\pscircle*(20,6){0.06}
\pscircle*(21,7){0.06}\pscircle*(22,8){0.06}\pscircle*(23,7){0.06}
\pscircle*(24,6){0.06}\pscircle*(25,6){0.06}\pscircle*(26,6){0.06}
\pscircle*(27,6){0.06}\pscircle*(28,7){0.06}\pscircle*(29,6){0.06}

\pscircle*(15,9){0.06}\pscircle*(16,10){0.06}\pscircle*(17,11){0.06}
\pscircle*(18,11){0.06}\pscircle*(19,10){0.06}\pscircle*(20,9){0.06}
\pscircle*(21,10){0.06}\pscircle*(22,11){0.06}\pscircle*(23,10){0.06}
\pscircle*(24,9){0.06}\pscircle*(25,9){0.06}\pscircle*(26,9){0.06}
\pscircle*(27,9){0.06}\pscircle*(28,10){0.06}\pscircle*(29,9){0.06}

\put(0.7,.25){$a$}\put(1.2,.25){$a$}\put(1.7,.25){$a$}\put(2.15,.05){$b$}
\put(3.3,.25){$a$}\put(3.8,.25){$a$}\put(4.2,.05){$b$}
\put(4.9,.1){$a$}\put(5.45,.1){$a$}\put(6,.1){$a$}\put(6.9,.05){$b$}
\put(4.9,.6){$w=a^8b^3$}

\put(1.1,2.15){$b$}\put(1.7,1.85){$a$}\put(2.1,1.65){$b$}
\put(3.3,1.85){$a$}\put(3.8,1.85){$a$}\put(4.2,1.65){$b$}
\put(4.9,1.7){$a$}\put(5.45,1.7){$a$}\put(6,1.7){$a$}\put(6.9,1.65){$b$}
\put(4.9,2.2){$w=a^6b^4$}

\put(0.7,3.45){$a$}\put(1.6,3.75){$b$}\put(2.15,3.25){$b$}
\put(3.3,3.45){$a$}\put(3.8,3.45){$a$}\put(4.2,3.25){$b$}
\put(4.9,3.3){$a$}\put(5.45,3.3){$a$}\put(6,3.3){$a$}\put(6.9,3.25){$b$}
\put(4.9,3.8){$w=a^6b^4$}

\put(1.2,5.55){$a$}\put(1.6,5.3){$b$}\put(2.15,4.8){$b$}
\put(3.3,5.05){$a$}\put(3.8,5.05){$a$}\put(4.2,4.8){$b$}
\put(4.9,4.85){$a$}\put(5.45,4.85){$a$}\put(6,4.85){$a$}\put(6.9,4.8){$b$}
\put(4.9,5.4){$w=a^6b^4$}

\put(8.6,.25){$a$}\put(9.1,.25){$a$}\put(9.65,.25){$a$}\put(10.05,.05){$b$}
\put(11.7,.5){$b$}\put(12.1,.05){$b$}
\put(12.8,.1){$a$}\put(13.35,.1){$a$}\put(13.9,.1){$a$}\put(14.8,.05){$b$}
\put(12.7,.6){$w=a^6b^4$}

\put(9,2.15){$b$}\put(9.65,1.85){$a$}\put(10.05,1.65){$b$}
\put(11.7,2.1){$b$}\put(12.1,1.65){$b$}
\put(12.8,1.7){$a$}\put(13.35,1.7){$a$}\put(13.9,1.7){$a$}\put(14.8,1.65){$b$}
\put(12.7,2.2){$w=a^4b^5$}

\put(8.6,3.45){$a$}\put(9.55,3.7){$b$}\put(10.05,3.2){$b$}
\put(11.7,3.7){$b$}\put(12.1,3.25){$b$}
\put(12.8,3.3){$a$}\put(13.35,3.3){$a$}\put(13.9,3.3){$a$}\put(14.8,3.25){$b$}
\put(12.7,3.8){$w=a^4b^5$}

\put(9.1,5.55){$a$}\put(9.55,5.3){$b$}\put(10.05,4.8){$b$}
\put(11.7,5.25){$b$}\put(12.1,4.8){$b$}
\put(12.8,4.85){$a$}\put(13.35,4.85){$a$}\put(13.9,4.85){$a$}\put(14.8,4.8){$b$}
\put(12.7,5.4){$w=a^4b^5$}

\psline*(8,14)(9,14.3)(11,14)
\psline*(13,14)(14,14.3)(15,14)

\put(4.8,7.7){$\mathbf{Q}_3$}\put(5.8,6.9){$b$}\put(7.2,7.7){$\mathbf{Q}_2$}\put(7.9,6.9){$b$}
\put(8.6,7){$a$}\put(9.1,7){$a$}\put(9.65,7){$a$}\put(10.55,6.9){$b$}
\put(5,6.4){$\phi(\mathbf{P})=\mathbf{u}\mathbf{Q}_3\mathbf{d}\mathbf{u}\mathbf{Q}_2\mathbf{d}\mathbf{hhh}\mathbf{u}\mathbf{d}$}

\end{pspicture}
\end{center}

\caption{\small An example of the bijection $\phi$ described in the proof of Corollary 3.2,
where $\mathbf{Q}_3$ has four possible cases and $\mathbf{Q}_2$ has two possible cases, that is $\mathbf{Q}_3=\mathbf{hhh}, \mathbf{udh}, \mathbf{hud}$ or $\mathbf{uhd}$ and $\mathbf{Q}_2=\mathbf{hh}$ or $\mathbf{ud}.$ }

\end{figure}
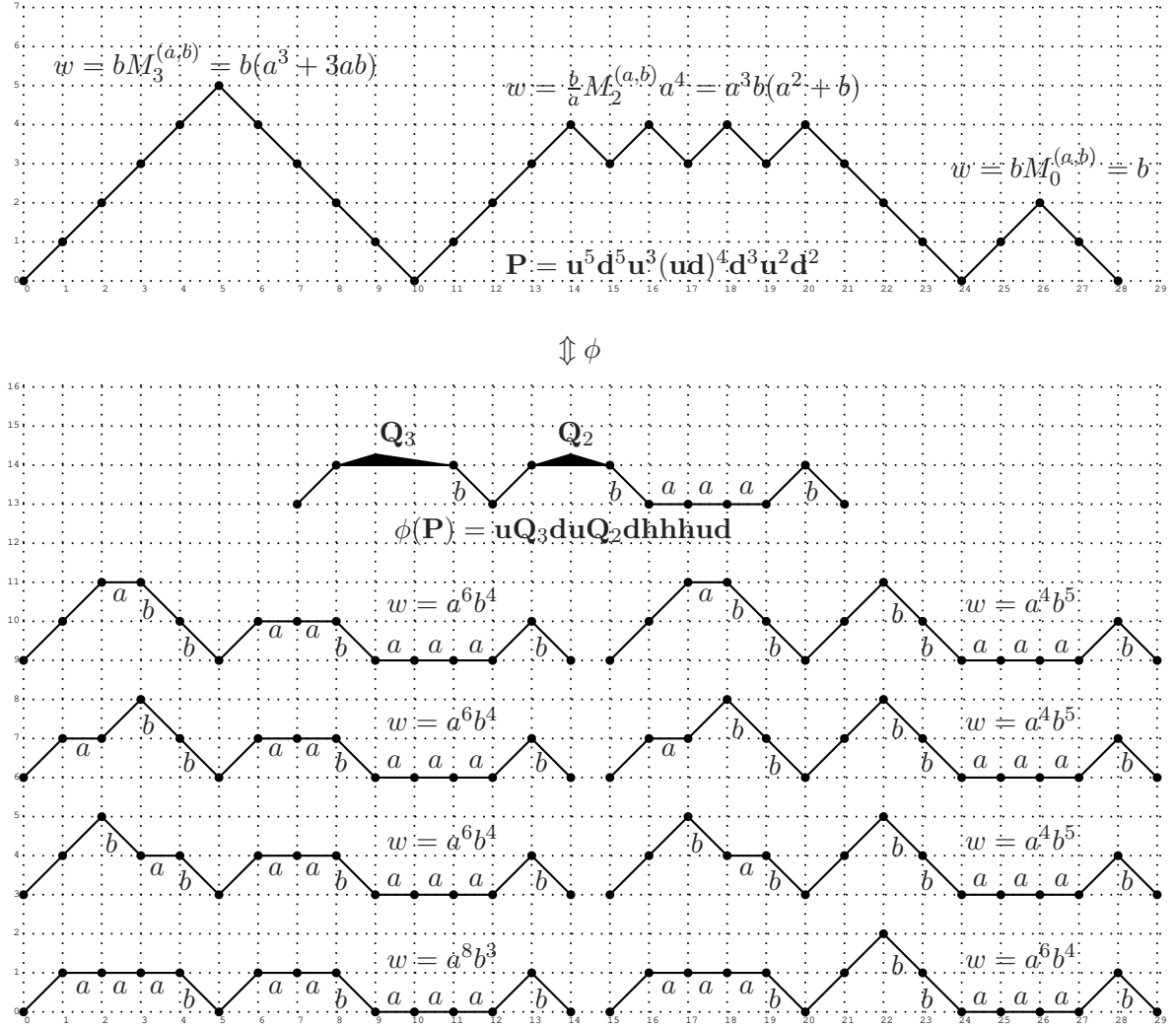

\subsection{The special cases related to Schr\"{o}der paths}

Let
$$R(x,q)=\sum_{n\geq 0}R_n(q)x^n=\frac{1-qx-\sqrt{(1-qx)^2-4x}}{2x}$$
be the generating function for the $q$-large Schr\"{o}der numbers $R_n(q)$, which counts the set $\mathcal{R}_n$ of $q$-large Schr\"{o}der paths of length $2n$ from $(0, 0)$ to $(2n, 0)$ in the first quadrant of the xy-plane that
consists of up steps $\mathbf{u}=(1, 1)$ with weight $1$, down steps $\mathbf{d}=(1, -1)$ with weight $1$ and horizontal steps $\mathbf{H}=(2, 0)$ with weight $q$ \cite{Yang}. When $q=0$, we obtain the Dyck paths, and $q=1$ we obtain the large Schr\"{o}der paths. Note that $R(x,q)$ has the relations
$$R(x,q)=1+qxR(x,q)+xR(x,q)^2=\frac{1}{1-qx-xR(x,q)}. $$

\begin{example}  The special case in Corollary \ref{coro2.2} when
$$\alpha(x)=(q+1)x, \beta(x)=\frac{R(x,q)-1}{q+1}, \gamma(x)=x(R(x,q)-1) $$
leads to
\begin{eqnarray*}
V_{\alpha, \beta}(x)=\frac{1-(q+1)x}{1-(q+1)x-x(R(x,q)-1)}=(1-(q+1)x)R(x,q),
\end{eqnarray*}
or equivalently, $V_n=R_n(q)-(q+1)R_{n-1}(q)$ for $n\geq 0$.
\end{example}

Let $\mathcal{Y}_n$ denote the set of $q$-large Schr\"{o}der paths of length $2n$ such that the first step is not an $\mathbf{H}$-step or the first two steps are not $\mathbf{ud}$. Then it is clear that $\mathcal{Y}_n$ is counted by $R_n(q)-(q+1)R_{n-1}(q)$. This shows the following corollary.
\begin{corollary}
There exists a simple bijection between the set $\mathcal{Y}_n$ and the set $\mathcal{V}_n$ with the weight functions $\alpha(x)=(q+1)x, \beta(x)=\frac{R(x,q)-1}{q+1}$ and $\gamma(x)=x(R(x,q)-1)$.
\end{corollary}
\pf Similar to the proof of Corollary 3.2, subject to the required weighted functions $\alpha(x)$ and $\gamma(x)$, for any $\mathbf{P}\in \mathcal{V}_n$ with $n\geq 2$, each primitive part of $\mathbf{P}$ is of the form
$\mathbf{u}^{k}(\mathbf{ud})^r\mathbf{d}^k$ for certain $k, r\geq 1$. Set $\mathbf{P}=\mathbf{P}_1\mathbf{P}_2$, where $\mathbf{P}_1$ is the first primitive part of length at least $4$ and  $\mathbf{P}_1=\mathbf{u}^{k}(\mathbf{ud})^r\mathbf{d}^k$ for certain $k, r\geq 1$.

When $r=1$, $\mathbf{P}_1=\mathbf{u}^{k+1}\mathbf{d}^{k+1}$ has weight $\gamma_{k}=R_{k}(q)$. Equivalently, we associate a $q$-large Schr\"{o}der path $\mathbf{u}\mathbf{Q}_{k}\mathbf{d}$ to $\mathbf{P}_1$, where $\mathbf{Q}_{k}\in \mathcal{R}_{k}$. In other words, we assign any $q$-large Schr\"{o}der path $\mathbf{u}\mathbf{Q}_{k}\mathbf{d}$ to a maximal pyramid of height $k+1$ at altitude $0$.

When $r\geq 2$, $\mathbf{P}_1$ has weight $\alpha_1^{r}\beta_{k}=(q+1)^{r-1}R_{k}(q)$. Equivalently, we associate a $q$-large Schr\"{o}der path $\mathbf{u}\mathbf{Q}_{k}\mathbf{d}\mathbf{s}_{1}\dots\mathbf{s}_{r-1}$ to $\mathbf{P}_1$, where $\mathbf{Q}_{k}\in \mathcal{R}_{k}$ and $\mathbf{s}_i=\mathbf{H}$ with weight $q$ or $\mathbf{s}_i=\mathbf{ud}$ with weight $1$ for $1\leq i\leq r-1$, if one notices that for $2\leq j\leq r$ the $j$-th peak $\mathbf{ud}$ in $\mathbf{P}_1$ with weight $q+1$ can be regarded as a peak $\mathbf{ud}$ with weight $q$ or $1$, and the peak $\mathbf{ud}$ with weight $q$ or $1$ corresponds respectively to an $\mathbf{H}$-step with weight $q$ or a peak $\mathbf{ud}$ with weight $1$ at $x$-axis.

Now we can recursively establish a bijection $\theta$ between $\mathcal{V}_n$ and $\mathcal{Y}_n$ as follows. For any $\mathbf{P}=\mathbf{P}_1\mathbf{P}_2\in \mathcal{V}_n$, where $\mathbf{P}_1=\mathbf{u}^{k}(\mathbf{ud})^r\mathbf{d}^k$ for certain $k, r\geq 1$, if $\mathbf{P}_1$ is associated with a $q$-large Schr\"{o}der path $\mathbf{u}\mathbf{Q}_{k}\mathbf{d}\mathbf{s}_{1}\dots\mathbf{s}_{r-1}$ as above, then $\theta(\mathbf{P})=\mathbf{u}\mathbf{Q}_{k}\mathbf{d}\mathbf{s}_{1}\dots\mathbf{s}_{r-1}\theta(\mathbf{P}_2)\in \mathcal{Y}_n$.

Conversely, for any $q$-large Schr\"{o}der path $\mathbf{Y}\in\mathcal{Y}_n$, $\mathbf{Y}$ can be written uniquely as $\mathbf{Y}=\mathbf{u}\mathbf{Q}_{k}\mathbf{d}\mathbf{s}_{1}\dots\mathbf{s}_{r-1}\mathbf{Y}_2$ for certain $k, r\geq 1$, where $\mathbf{Q}_{k}\in \mathcal{R}_{k}$ and $\mathbf{Y}_2\in \mathcal{Y}_{n-k-r}$, the inverse $\theta^{-1}$ is built by $\theta^{-1}(\mathbf{Y})=\mathbf{P}_1\theta^{-1}(\mathbf{Y}_2)\in \mathcal{V}_n$, where $\mathbf{P}_1=\mathbf{u}^{k}(\mathbf{ud})^r\mathbf{d}^k$ is associated with a $q$-large Schr\"{o}der path $\mathbf{u}\mathbf{Q}_{k}\mathbf{d}\mathbf{s}_{1}\dots\mathbf{s}_{r-1}$ and such paths have the total weight $(q+1)^{r-1}R_{k}(q)$.    \qed  \vskip0.2cm

In order to give a more intuitive view on the bijection $\theta$, we present a pictorial description of
$\theta$ for the case $\mathbf{P}=\mathbf{u}^3\mathbf{d}^3\mathbf{u}(\mathbf{u}\mathbf{d})^3\mathbf{d}$ with weight $R_{2}(q)\cdot\frac{1}{q+1}R_{1}(q)(q+1)^{3}=(q+2)(q+1)^4$. Clearly, $\theta(\mathbf{P})=\mathbf{u}\mathcal{R}_{2}\mathbf{d}\mathbf{u}\mathcal{R}_{1}\mathbf{d}\mathbf{s}_1\mathbf{s}_2$ with $\mathbf{s}_i\in \{\mathbf{H}, \mathbf{ud}\}$. Actually, if the first primitive part $\mathbf{u}^3\mathbf{d}^3$ of $\mathbf{P}$ is associated with a primitive $q$-large Schr\"{o}der path $\mathbf{u}\mathbf{Q}_{2}\mathbf{d}$, the beginning $\mathbf{u}$-segment $\mathbf{u}^2$ of the second primitive part $\mathbf{u}(\mathbf{u}\mathbf{d})^3\mathbf{d}$ of $\mathbf{P}$ is associated with a primitive $q$-large Schr\"{o}der path $\mathbf{u}\mathbf{Q}_{1}\mathbf{d}$, where $\mathbf{Q}_i\in \mathcal{R}_{i}$ for $i=1, 2$, and
the last two $\mathbf{ud}$-peaks in $\mathbf{u}(\mathbf{u}\mathbf{d})^3\mathbf{d}$ are weighted respectively by $q$ and $1$, then $\theta(\mathbf{P})=\mathbf{u}\mathbf{Q}_2\mathbf{d}\mathbf{u}\mathbf{Q}_1\mathbf{d}\mathbf{Hud}\in \mathcal{Y}_{14}$. See Figure 3.

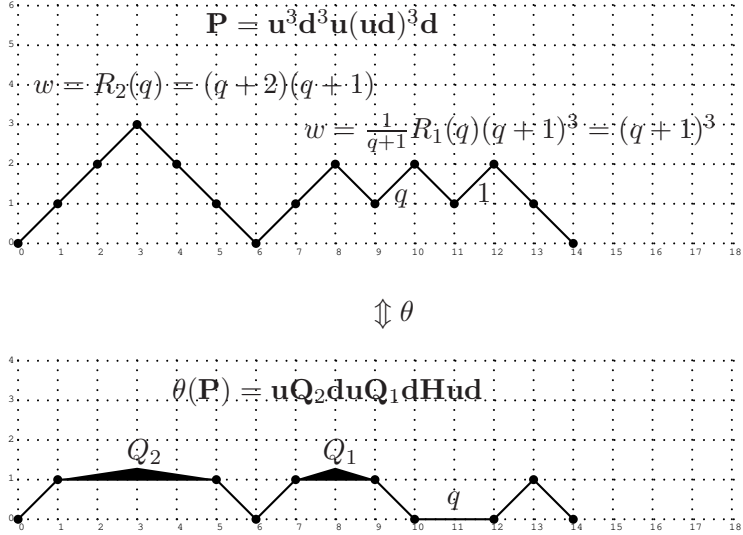
\begin{figure}[h] \setlength{\unitlength}{0.5mm}

\begin{center}
\begin{pspicture}(10,3.5)
\psset{xunit=15pt,yunit=15pt}\psgrid[subgriddiv=1,griddots=4,
gridlabels=4pt](0,0)(18,6)

\psline(0,0)(3,3)(6,0)(8,2)(9,1)(10,2)(11,1)(12,2)(14,0)

\pscircle*(0,0){0.06}\pscircle*(1,1){0.06}\pscircle*(2,2){0.06}
\pscircle*(3,3){0.06}\pscircle*(4,2){0.06}\pscircle*(5,1){0.06}
\pscircle*(6,0){0.06}\pscircle*(7,1){0.06}\pscircle*(8,2){0.06}
\pscircle*(9,1){0.06}\pscircle*(10,2){0.06}\pscircle*(11,1){0.06}
\pscircle*(12,2){0.06}\pscircle*(13,1){0.06}\pscircle*(14,0){0.06}

\put(5,0.55){$q$} \put(6.1,0.55){$1$}

\put(0.2,2){$\tiny w=R_{2}(q)=(q+2)(q+1)$} \put(3.8,1.4){$w=\frac{1}{q+1}R_{1}(q)(q+1)^{3}=(q+1)^3$}
\put(2.5,2.8){$\mathbf{P}=\mathbf{u}^3\mathbf{d}^3\mathbf{u}(\mathbf{u}\mathbf{d})^3\mathbf{d}$}

\end{pspicture}
\end{center}\vskip0.5cm

$\Updownarrow \theta$
\vskip.5cm

\begin{center}
\begin{pspicture}(10,2)
\psset{xunit=15pt,yunit=15pt}\psgrid[subgriddiv=1,griddots=4,
gridlabels=4pt](0,0)(18,4)

\psline(0,0)(1,1)(4,1)(5,1)(6,0)(7,1)(9,1)(10, 0)(12,0)(13,1)(14,0)

\psline*(1,1)(3,1.3)(5,1)

\psline*(7,1)(8,1.3)(9,1)

\pscircle*(0,0){0.06}\pscircle*(1,1){0.06}\pscircle*(5,1){0.06}
\pscircle*(6,0){0.06}\pscircle*(7,1){0.06}
\pscircle*(9,1){0.06}\pscircle*(10,0){0.06}
\pscircle*(12,0){0.06}\pscircle*(13,1){0.06}\pscircle*(14,0){0.06}

\put(1.45,0.8){$Q_2$}\put(4.05,0.8){$Q_1$}\put(5.7,0.2){$q$}

\put(2.05,1.6){$\theta(\mathbf{P})=\mathbf{u}\mathbf{Q}_2\mathbf{d}\mathbf{u}\mathbf{Q}_1\mathbf{d}\mathbf{Hud}$}

\end{pspicture}
\end{center}

\caption{\small An example of the bijection $\theta$ described in the proof of Corollary 3.4,
where $\mathbf{Q}_2$ has six possible cases and $\mathbf{Q}_1$ has two possible cases, that is $\mathbf{Q}_2=\mathbf{HH}, \mathbf{Hud}, \mathbf{uHd}$, $\mathbf{udH}$, $\mathbf{udud}$ or $\mathbf{uudd}$ and $\mathbf{Q}_1=\mathbf{H}$ or $\mathbf{ud}.$ }

\end{figure}

There is another special case related to Schr\"{o}der paths. Let
$$S(x,q)=\sum_{n\geq 0}S_n(q)x^n=\frac{1+qx-\sqrt{(1+qx)^2-4(1+q)x}}{2(1+q)x}$$
be the generating function for the $q$-small Schr\"{o}der numbers $S_n(q)$, which count the set $\mathcal{S}_n$ of $q$-small Schr\"{o}der paths, i.e., $q$-large Schr\"{o}der paths of length $2n$ from $(0, 0)$ to $(2n, 0)$ without horizontal steps at $x$-axis \cite{Yang}. When $q=0$, we obtain the Dyck paths, and when $q=1$ we obtain the small Schr\"{o}der paths. Note that $S(x,q)$ has the relations
$$S(x,q)=1-qxS(x,q)+(1+q)xS(x,q)^2=\frac{1}{1+qx-(1+q)xS(x,q)}. $$

\begin{example}  The special case in Corollary \ref{coro2.2} when
$$\alpha(x)=x, \beta(x)=(q+1)(S(x,q)-1), \gamma(x)=(q+1)x(S(x,q)-1) $$
leads to
\begin{eqnarray*}
V_{\alpha, \beta}(x)=\frac{1-x}{1-x-(q+1)x(S(x,q)-1)}=(1-x)S(x,q),
\end{eqnarray*}
or equivalently, $V_n=S_n(q)-S_{n-1}(q)$ for $n\geq 0$.
\end{example}

Let $\mathcal{Z}_n$ denote the set of $q$-small Schr\"{o}der paths of length $2n$ such that the first two steps are not $\mathbf{ud}$. Then it is clear that $\mathcal{Z}_n$ is counted by $S_n(q)-S_{n-1}(q)$. This shows the following corollary.
\begin{corollary}
There exists a simple bijection between the set $\mathcal{Z}_n$ and the set $\mathcal{V}_n$ with the weight functions $\alpha(x)=x, \beta(x)=(q+1)(S(x,q)-1)$ and $\gamma(x)=(q+1)x(S(x,q)-1)$.
\end{corollary}
\pf The proof is similar to that of Corollary 3.4, which is left to the interested readers. \qed \vskip0.2cm

\subsection{The special cases related to Narayana polynomials}
Let $N(x,t)=\sum_{n\geq 0}N_n(t)x^n$ be the generating function for the Narayana polynomials
$$N_n(t)=\sum_{i=1}^{n}\frac{1}{n}\binom{n}{i}\binom{n}{i-1}t^{i}$$
with $N_0(t)=1$ and $N_1(t)=t$. The Narayana polynomials $N_n(t)$ count the weighted Dyck paths $\mathcal{D}^{*}_n$ of length $2n$ such that each peak is weighted by $t$ and other steps are weighted by $1$ \cite{Deutsch99}. Both the recurrences and the explicit formula for $N(x,t)$ are given respectively by
\begin{eqnarray*}
N(x,t) \hskip-.22cm &=& \hskip-.22cm 1+(t-1)xN(x,t)+xN(x,t)^2=\frac{1}{1-(t-1)x-xN(x,t)} \\
       \hskip-.22cm &=& \hskip-.22cm  \frac{1-(t-1)x-\sqrt{1-2(1+t)x+(1-t)^2x^2}}{2x}.
\end{eqnarray*}
There is a closely relation between the Narayana polynomials $N_n(t)$ and the $q$-large Schr\"{o}der numbers $R_n(q)$, in fact, $R_n(q)=N_n(q+1)$ \cite{Bonin}.

\begin{example}
The special case in Corollary \ref{coro2.2} when
$$\alpha(x)=tx, \beta(x)=\frac{N(x,t)-1}{t}, \gamma(x)=x(N(x,t)-1) $$
leads to
\begin{eqnarray*}
V_{\alpha, \beta}(x)=\frac{1-tx}{1-tx-x(N(x,t)-1)}=(1-tx)N(x,t),
\end{eqnarray*}
or equivalently, $V_n=N_n(t)-tN_{n-1}(t)$ for $n\geq 0$.

Note that $N_n(1)=C_{n}$ and $N_n(q+1)=R_n(q)$, the results above can further generate $V_n=C_{n}-C_{n-1}$ in $t=1$ and $V_n=R_n(q)-(q+1)R_{n-1}(q)$ in $t=q+1$ which is the result in Example 3.3 subject to the $q$-large Schr\"{o}der paths.

\end{example}

Let $\mathcal{U}_n$ denote the subset of $\mathcal{D}^{*}_n$ such that the first two steps are not $\mathbf{ud}$. Clearly, $\mathcal{U}_n$ is counted by $N_n(t)-tN_{n-1}(t)$. This signifies the following corollary.
\begin{corollary}
There exists a simple bijection between the set $\mathcal{U}_n$ and the set $\mathcal{V}_n$ with the weight functions $\alpha(x)=tx, \beta(x)=\frac{N(x,t)-1}{t}, \gamma(x)=x(N(x,t)-1)$.
\end{corollary}

\pf Similar to the proof of Corollary 3.2, for any $\mathbf{P}\in \mathcal{V}_n$ with $n\geq 2$, $\mathbf{P}$ can be written uniquely as $\mathbf{P}=\mathbf{P}_1\mathbf{P}_2$, where $\mathbf{P}_1$ are the first primitive part of length at least $4$ and $\mathbf{P}_1=\mathbf{u}^{k}(\mathbf{ud})^r\mathbf{d}^k$ for certain $k, r\geq 1$.

When $r=1$, $\mathbf{P}_1=\mathbf{u}^{k+1}\mathbf{d}^{k+1}$ has weight $\gamma_{k}=N_{k}(t)$. Equivalently, we associate a primitive Dyck path $\mathbf{u}\mathbf{Q}_{k}\mathbf{d}$ to $\mathbf{P}_1$, where $\mathbf{Q}_{k}\in \mathcal{D}^{*}_k$. In other words, we assign any primitive Dyck path $\mathbf{u}\mathbf{Q}_{k}\mathbf{d}\in \mathcal{D}^{*}_{k+1}$ to a maximal pyramid of height $k+1$ at altitude $0$.

When $r\geq 2$, $\mathbf{P}_1$ has weight $\alpha_1^{r}\beta_{k}=t^{r-1}N_{k}(t)$. Equivalently, we associate a Dyck path $\mathbf{u}\mathbf{Q}_{k}\mathbf{d}(\mathbf{ud})^{r-1}\in \mathcal{D}^{*}_{k+r}$ to $\mathbf{P}_1$.

Now we can recursively establish a bijection $\rho$ between $\mathcal{V}_n$ and $\mathcal{U}_n$ as follows. For any $\mathbf{P}=\mathbf{P}_1\mathbf{P}_2\in \mathcal{V}_n$, where $\mathbf{P}_1=\mathbf{u}^{k}(\mathbf{ud})^r\mathbf{d}^k$ for certain $k, r\geq 1$, if $\mathbf{P}_1$ is associated by a Dyck path $\mathbf{u}\mathbf{Q}_{k}\mathbf{d}(\mathbf{ud})^{r-1}\in \mathcal{D}^{*}_{k+r}$, then $\rho(\mathbf{P})=\mathbf{u}\mathbf{Q}_{k}\mathbf{d}(\mathbf{ud})^{r-1}\rho(\mathbf{P}_2)\in \mathcal{U}_n$.

Conversely, for any Dyck path $\mathbf{U}\in\mathcal{U}_n$, $\mathbf{U}$ can be written uniquely as $\mathbf{U}=\mathbf{u}\mathbf{Q}_{k}\mathbf{d}(\mathbf{ud})^{r-1}\mathbf{U}_2$ for certain $k, r\geq 1$, where $\mathbf{u}\mathbf{Q}_{k}\mathbf{d}(\mathbf{ud})^{r-1}\in \mathcal{D}^{*}_{k+r}$ and $\mathbf{U}_2\in \mathcal{U}_{n-k-r}$, the inverse $\rho^{-1}$ is built by $\rho^{-1}(\mathbf{U})=\mathbf{P}_1\rho^{-1}(\mathbf{U}_2)\in \mathcal{V}_n$, where $\mathbf{P}_1=\mathbf{u}^{k}(\mathbf{ud})^r\mathbf{d}^k$ is associated with a weighted Dyck path $\mathbf{u}\mathbf{Q}_{k}\mathbf{d}(\mathbf{ud})^{r-1}\in \mathcal{D}^{*}_{k+r}$ and such weighted Dyck paths have the total weight $t^{r-1}N_{k}(t)$.    \qed  \vskip0.2cm

In order to give a more intuitive view on the bijection $\rho$, we present a pictorial description of
$\rho$ for the case $\mathbf{P}=\mathbf{u}^3\mathbf{d}^3\mathbf{u}^2(\mathbf{u}\mathbf{d})^4\mathbf{d}^2$ with weight $N_{2}(t)N_{2}(t)t^3=(t+t^2)^2t^3$. Clearly, $\rho(\mathbf{P})=\mathbf{u}\mathcal{D}^{*}_{2}\mathbf{d}\mathbf{u}\mathcal{D}^{*}_{2}\mathbf{d}(\mathbf{ud})^3$ with $\mathbf{ud}$ weighted by $t$. Precisely, if the first primitive part $\mathbf{u}^3\mathbf{d}^3$ of $\mathbf{P}$ is associated with a primitive Dyck path $\mathbf{u}\mathbf{Q}_{2}\mathbf{d}$, the beginning $\mathbf{u}$-segment $\mathbf{u}^3$ of the second primitive part $\mathbf{u}^2(\mathbf{u}\mathbf{d})^4\mathbf{d}^2$ of $\mathbf{P}$ is associated with a primitive Dyck path $\mathbf{u}\mathbf{Q}'_{2}\mathbf{d}$, where $\mathbf{Q}_2, \mathbf{Q}'_2\in \mathcal{D}^{*}_{2}$, and each of
the last three $\mathbf{ud}$-peaks in $\mathbf{u}^2(\mathbf{u}\mathbf{d})^4\mathbf{d}^2$ is weighted by $t$, then $\rho(\mathbf{P})=\mathbf{u}\mathbf{Q}_2\mathbf{d}\mathbf{u}\mathbf{Q}'_2\mathbf{d}\mathbf{ududud}\in \mathcal{U}_{18}$. See Figure 4.

\begin{figure}[h] \setlength{\unitlength}{0.5mm}

\begin{center}
\begin{pspicture}(10,3.5)
\psset{xunit=15pt,yunit=15pt}\psgrid[subgriddiv=1,griddots=4,
gridlabels=4pt](0,0)(18,6)

\psline(0,0)(3,3)(6,0)(8,2)(9,3)(10,2)(11,3)(12,2)(13,3)(14,2)(15,3)(18,0)

\pscircle*(0,0){0.06}\pscircle*(1,1){0.06}\pscircle*(2,2){0.06}
\pscircle*(3,3){0.06}\pscircle*(4,2){0.06}\pscircle*(5,1){0.06}
\pscircle*(6,0){0.06}\pscircle*(7,1){0.06}\pscircle*(8,2){0.06}
\pscircle*(9,3){0.06}\pscircle*(10,2){0.06}\pscircle*(11,3){0.06}
\pscircle*(12,2){0.06}\pscircle*(13,3){0.06}\pscircle*(14,2){0.06}
\pscircle*(15,3){0.06}\pscircle*(16,2){0.06}\pscircle*(17,1){0.06}
\pscircle*(18,0){0.06}

\put(4.6,1.1){$t$} \put(5.7,1.1){$t$}\put(6.8,1.1){$t$}\put(7.9,1.1){$t$}

\put(0.2,2){$\tiny w=N_{2}(t)=(t+t^2)$} \put(4.4,2){$w=\frac{1}{t}N_{2}(t)t^{4}=(t+t^2)t^3$}
\put(3,2.8){$\mathbf{P}=\mathbf{u}^3\mathbf{d}^3\mathbf{u^2}(\mathbf{u}\mathbf{d})^4\mathbf{d}^2$}

\end{pspicture}
\end{center}\vskip0.5cm

$\Updownarrow \rho$
\vskip.5cm

\begin{center}
\begin{pspicture}(10,2)
\psset{xunit=15pt,yunit=15pt}\psgrid[subgriddiv=1,griddots=4,
gridlabels=4pt](0,0)(18,4)

\psline(0,0)(1,1)(4,1)(5,1)(6,0)(7,1)(9,1)(10,1)(11,1)(12,0)(13,1)(14,0)(15,1)(16,0)(17,1)(18,0)

\psline*(1,1)(3,1.3)(5,1)

\psline*(7,1)(9,1.3)(11,1)

\pscircle*(0,0){0.06}\pscircle*(1,1){0.06}\pscircle*(5,1){0.06}
\pscircle*(6,0){0.06}\pscircle*(7,1){0.06}
\pscircle*(11,1){0.06}\pscircle*(12,0){0.06}\pscircle*(13,1){0.06}
\pscircle*(14,0){0.06}\pscircle*(15,1){0.06}\pscircle*(16,0){0.06}
\pscircle*(17,1){0.06}\pscircle*(18,0){0.06}

\put(1.45,0.8){$Q_2$}\put(4.55,0.8){$Q'_2$}\put(6.75,0.1){$t$}\put(7.85,0.1){$t$}\put(8.95,0.1){$t$}

\put(2.05,1.6){$\rho(\mathbf{P})=\mathbf{u}\mathbf{Q}_2\mathbf{d}\mathbf{u}\mathbf{Q}'_2\mathbf{d}\mathbf{ududud}$}

\end{pspicture}
\end{center}

\caption{\small An example of the bijection $\rho$ described in the proof of Corollary 3.8,
where each of $\mathbf{Q}_2$ and $\mathbf{Q}'_2$ has two possible cases, that is  $\mathbf{Q}_2, \mathbf{Q}'_2\in \{\mathbf{uudd}, \mathbf{udud}\}$. }

\end{figure}
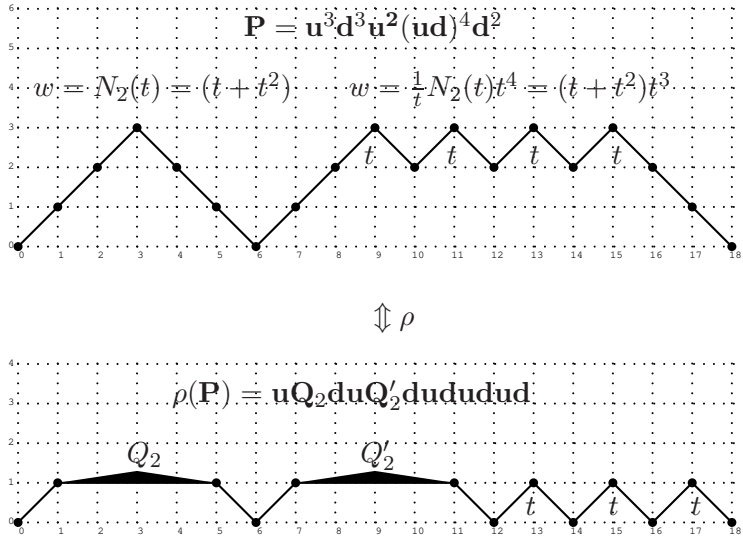

Now we can consider another case related to Narayana polynomials.
\begin{example} Let $f(x,t)=\sum_{n\geq 0}\frac{1}{t}N_{n+1}(t)x^n=\frac{N(x,t)-1}{tx}$, the recurrence for $f(x,t)$ is
\begin{eqnarray*}
f(x,t) \hskip-.22cm &=& \hskip-.22cm (1+xf(x,t))(1+xtf(x,t))=\frac{1}{1-(1+t)x-x^2tf(x,t)}.
\end{eqnarray*}
Then, the special case in Corollary \ref{coro2.2} when
$$\alpha(x)=(1+t)x, \beta(x)=\frac{tx}{1+t}f(x,t), \gamma(x)=tx^2f(x,t) $$
leads to
\begin{eqnarray*}
V_{\alpha, \beta}(x)=\frac{1-(1+t)x}{1-(1+t)x-x^2tf(x,t)}=(1-(1+t)x)f(x,t),
\end{eqnarray*}
or equivalently, $V_n=\frac{1}{t}(N_{n+1}(t)-(t+1)N_{n}(t))$ for $n\geq 0$.

The cases $t=1$ and $t=q+1$ give $V_n=C_{n+1}-2C_{n}$ and $V_n=\frac{1}{q+1}(R_{n+1}(q)-(q+2)R_{n}(q))$ respectively.

\end{example}

Let $\mathcal{\bar{U}}_n$ denote the set of Dyck paths of length $2n$ with peaks at level $1$ marked by $t+1$ and peaks at level $\geq 2$ marked by $t$ such that the first two steps are not $\mathbf{ud}$. Then, we have
\begin{corollary}
There exists a simple bijection $\psi$ between the set $\mathcal{\bar{U}}_n$ and the set $\mathcal{V}_n$ with the weight functions $\alpha(x)=(t+1)x, \beta(x)=\frac{tx}{1+t}f(x,t), \gamma(x)=tx^2f(x,t)$. Hence, $\mathcal{\bar{U}}_n$ is counted by $\frac{1}{t}(N_{n+1}(t)-(t+1)N_{n}(t))$.
\end{corollary}

\pf Similar to the proof of Corollary 3.2, for any $P\in \mathcal{V}_n$ with $n\geq 2$, $\mathbf{P}$ can be written uniquely as $\mathbf{P}=\mathbf{P}_1\mathbf{P}_2$, where $\mathbf{P}_1$ are primitive and of length at least $2$. In fact, each $\mathbf{P}_1=\mathbf{u}^{k}(\mathbf{ud})^r\mathbf{d}^k$ for certain $k, r\geq 1$.

When $r=1$, $\mathbf{P}_1=\mathbf{u}^{k+1}\mathbf{d}^{k+1}$ has weight $\gamma_{k}=N_{k}(t)$. Equivalently, we associate a Dyck path $\mathbf{u}\mathbf{Q}_{k-1}\mathbf{d}$ to $\mathbf{P}_1$, where $\mathbf{Q}_{k}\in \mathcal{D}^{*}_{k}$. In other words, we assign any primitive Dyck path $\mathbf{u}\mathbf{Q}_{k}\mathbf{d}\in \mathcal{D}^{*}_{k+1}$ to a maximal pyramid of height $k+1$ at altitude $0$.

When $r\geq 2$, $P_1$ has weight $\alpha_1^{r}\beta_{k}=(t+1)^{r-1}N_{k}(t)$. Equivalently, we associate a Dyck path $\mathbf{u}\mathbf{Q}_{k}\mathbf{d}(\mathbf{ud})^{r-1}$ with each peak at level $1$ weighted by $t+1$ to $\mathbf{P}_1$, where $\mathbf{Q}_{k}\in \mathcal{D}^{*}_{k}$.

Now we can recursively establish a bijection $\psi$ between $\mathcal{V}_n$ and $\mathcal{\bar{U}}_n$ as follows. For any $\mathbf{P}=\mathbf{P}_1\mathbf{P}_2\in \mathcal{V}_n$, where $\mathbf{P}_1=\mathbf{u}^{k}(\mathbf{ud})^r\mathbf{d}^k$ for certain $k, r\geq 1$, if $\mathbf{P}_1$ is associated by a weighted Dyck path $\mathbf{u}\mathbf{Q}_{k}\mathbf{d}(\mathbf{ud})^{r-1}$ as above, then $\psi(\mathbf{P})=\mathbf{u}\mathbf{Q}_{k}\mathbf{d}(\mathbf{ud})^{r-1}\psi(\mathbf{P}_2)\in \mathcal{\bar{U}}_n$.

Conversely, for any Dyck path $\mathbf{U}\in\mathcal{\bar{U}}_n$, $\mathbf{U}$ can be written uniquely as $\mathbf{U}=\mathbf{u}\mathbf{Q}_{k}\mathbf{d}(\mathbf{ud})^{r-1}\mathbf{U}_2$ for certain $k, r\geq 1$, where $\mathbf{Q}_{k}\in \mathcal{D}^{*}_{k}$ and $\mathbf{U}_2\in \mathcal{U}_{n-k-r}$, the inverse $\psi^{-1}$ is built by $\psi^{-1}(\mathbf{U})=\mathbf{P}_1\psi^{-1}(\mathbf{U}_2)\in \mathcal{V}_n$, where $\mathbf{P}_1=\mathbf{u}^{k}(\mathbf{ud})^r\mathbf{d}^k$ associated with a Dyck path $\mathbf{u}\mathbf{Q}_{k}\mathbf{d}(\mathbf{ud})^{r-1}$ such that $\mathbf{Q}_{k}\in \mathcal{D}^{*}_{k}$ and each peak $\mathbf{ud}$ at level $1$ weighted by $t+1$, such Dyck paths have the total weight $(t+1)^{r-1}N_{k}(t)$.

Clearly, according to the bijection $\psi$, $\mathcal{\bar{U}}_n$ is counted by $\frac{1}{t}(N_{n+1}(t)-(t+1)N_{n}(t))$.  \qed  \vskip0.2cm

The pictorial description of $\psi$ is similar to that of $\rho$ which is left to interested readers.

\subsection{The special cases related to the Chebyshev polynomials} Let
$$U(x,t)=\sum_{n\geq 0}U_n(t)x^n=\frac{1}{1-2tx+x^2}$$
be the generating function for the Chebyshev polynomials of the second kind \cite[Page 49]{Comtet}.
For nonnegative numbers $a, b, c, d$ with $a> b$, let
$$\alpha(x)=\frac{(a-b)x}{1-bx}, \beta(x)=\frac{cx}{1-dx}, \gamma(x)=\frac{(a-b)cx^2}{(1-bx)(1-dx)}.$$
Then the special case in Corollary \ref{coro2.2} produces
\begin{eqnarray}\label{eqn 3.1}
V_{\alpha, \beta}(x)=1+\frac{(a-b)cx^2}{1-(a+d)x+(ad-(a-b)c)x^2}.
\end{eqnarray}

$\mathbf{Case\ 1}$. When $ad=(a-b)c$, we have $V_n=ad(a+d)^{n-2}$ for $n\geq 2$ and $V_0=1, V_1=0.$

$\mathbf{Case\ 2}$. When $ad=(a-b)c+1$, we have
\begin{eqnarray*}
V_{\alpha, \beta}(x)=1+\frac{(ad-1)x^2}{1-(a+d)x+x^2}=1+(ad-1)x^2\sum_{n=0}^{\infty}U_n\Big(\frac{a+d}{2}\Big)x^{n},
\end{eqnarray*}
and $V_n=(ad-1)U_{n-2}(\frac{a+d}{2})$ for $n\geq 2$ with $V_0=1, V_1=0$.

$\mathbf{Case\ 3}$. When $a+d=3$ and $ad=(a-b)c+1$, we have
\begin{eqnarray*}
V_{\alpha, \beta}(x)=1+\frac{(ad-1)x^2}{1-3x+x^2}=1+(ad-1)\sum_{n=1}^{\infty}F_{2n}x^{n+1},
\end{eqnarray*}
and $V_n=(ad-1)F_{2n-2}$ for $n\geq 1$ with $V_0=1$, where $F_n$ denote the Fibonacci numbers defined by $F_n=F_{n-1}+F_{n-2}$ with $F_0=1, F_1=1.$

There is another case related to the Chebyshev polynomials. For $a, b, c\geq 0$, let
$$\alpha(x)=\frac{2bx(1-ax)}{1-2cx+x^2}, \beta(x)=\frac{ax}{1-ax}, \gamma(x)=\frac{2abx^2}{1-2cx+x^2}.$$
Then the special case in Corollary \ref{coro2.2} produces
\begin{eqnarray*}
V_{\alpha, \beta}(x)=1+\frac{2abx^2}{1-2(b+c)x+x^2}=1+2abx^2\sum_{n=0}^{\infty}U_n(b+c)x^{n},
\end{eqnarray*}
and $V_n=2abU_{n-2}(b+c)$ for $n\geq 2$ with $V_0=1, V_1=0$.

\subsection{The special case related to Delannoy paths}

Let $D(x)=\sum_{n\geq 0}D_nx^n=\frac{1}{\sqrt{1-6x+x^2}}$ be the generating function for central Delannoy numbers
$$D_n=\sum_{i=0}^{n}\binom{n}{i}\binom{n+i}{i}=\sum_{i=0}^{n}\binom{n}{i}^{2}2^{i},$$
which counts the number of Delannoy paths from $(0,0)$ to $(2n,0)$ using up steps $\mathbf{u}=(1, 1)$, down steps $\mathbf{d}=(1, -1)$ and horizontal steps $\mathbf{H}=(2, 0)$ \cite{Banderier}. Note that $\sum_{i=0}^{n-1}D_iD_{n-i-1}$ counts the total number of $\mathbf{H}$-steps of Delannoy paths on the X-axis from $(0,0)$ to $(2n,0)$ for $n\geq 1$.

When $a+d=6$ and $ad=(a-b)c+1$ in (\ref{eqn 3.1}), we have
\begin{eqnarray*}
V_{\alpha, \beta}(x)=1+\frac{(ad-1)x^2}{1-6x+x^2},
\end{eqnarray*}
and $V_n=(ad-1)\sum_{i=0}^{n-2}D_iD_{n-2-i}$ for $n\geq 2$ with $V_0=1, V_1=0.$ It is interesting that $V_n$ in Table 3.1 are invariant under the following assignment for $a, b, c$ and $d$,
\begin{eqnarray*}
{\rm (I)} \left\{\begin{array}{ccc}
 a \hskip-.22cm &=&\hskip-.22cm 4  \\
 b \hskip-.22cm &=&\hskip-.22cm 3  \\
 c \hskip-.22cm &=&\hskip-.22cm 7  \\
 d \hskip-.22cm &=&\hskip-.22cm 2
\end{array}\right.
\mbox{and} \
\left\{\begin{array}{ccc}
 a \hskip-.22cm &=&\hskip-.22cm 2  \\
 b \hskip-.22cm &=&\hskip-.22cm 1  \\
 c \hskip-.22cm &=&\hskip-.22cm 7  \\
 d \hskip-.22cm &=&\hskip-.22cm 4
\end{array}\right.,\
{\rm (II)} \left\{\begin{array}{ccc}
 a \hskip-.22cm &=&\hskip-.22cm 5  \\
 b \hskip-.22cm &=&\hskip-.22cm 4  \\
 c \hskip-.22cm &=&\hskip-.22cm 4  \\
 d \hskip-.22cm &=&\hskip-.22cm 1
\end{array}\right.,
\left\{\begin{array}{ccc}
 a \hskip-.22cm &=&\hskip-.22cm 5  \\
 b \hskip-.22cm &=&\hskip-.22cm 1  \\
 c \hskip-.22cm &=&\hskip-.22cm 1  \\
 d \hskip-.22cm &=&\hskip-.22cm 1
\end{array}\right.
\mbox{and} \
\left\{\begin{array}{ccc}
 a \hskip-.22cm &=&\hskip-.22cm 1  \\
 b \hskip-.22cm &=&\hskip-.22cm 0  \\
 c \hskip-.22cm &=&\hskip-.22cm 4  \\
 d \hskip-.22cm &=&\hskip-.22cm 5
\end{array}\right., \\
{\rm (III)} \left\{\begin{array}{ccc}
 a \hskip-.22cm &=&\hskip-.22cm 3  \\
 b \hskip-.22cm &=&\hskip-.22cm 2  \\
 c \hskip-.22cm &=&\hskip-.22cm 8  \\
 d \hskip-.22cm &=&\hskip-.22cm 3
\end{array}\right.
\mbox{and} \
\left\{\begin{array}{ccc}
 a \hskip-.22cm &=&\hskip-.22cm 3  \\
 b \hskip-.22cm &=&\hskip-.22cm 1  \\
 c \hskip-.22cm &=&\hskip-.22cm 4  \\
 d \hskip-.22cm &=&\hskip-.22cm 3
\end{array}\right..
\end{eqnarray*}

\begin{center}
\small\begin{eqnarray*}
\begin{array}{|c|c|c|c|}\hline
                    & \rm Pair\ (I)                                & \rm Pair\ (II)                                & \rm  Pair\ (III)                              \\[3pt]\hline
V(x)                &\displaystyle 1+\frac{7x^2}{1-6x+x^2}         & \displaystyle 1+\frac{4x^2}{1-6x+x^2}         & \displaystyle 1+\frac{8x^2}{1-6x+x^2}          \\[5pt]\hline
V_n,(n\geq 2)       &\displaystyle 7\sum_{i=0}^{n-2}D_iD_{n-2-i}   & \displaystyle 4\sum_{i=0}^{n-2}D_iD_{n-2-i}   & \displaystyle  8\sum_{i=0}^{n-2}D_iD_{n-2-i}      \\[5pt] \hline
\end{array}
\end{eqnarray*}
Table 3.1. $V(x)$ and $V_{n}$ according to Pairs (I), (II) and (III).
\end{center}

\begin{corollary}
Let $\mathcal{V}_n^{(a,b,c,d)}$ be the set of $\mathcal{V}_n$ with the weight functions $\alpha(x)=\frac{(a-b)x}{1-bx}$, $\beta(x)=\frac{cx}{1-dx}$, $\gamma(x)=\frac{(a-b)cx^2}{(1-bx)(1-dx)}$.
Then there exist bijections between $\mathcal{V}_n^{(4,3,7,2)}$ and $\mathcal{V}_n^{(2,1,7,4)}$, between $\mathcal{V}_n^{(5,4,4,1)}$, $\mathcal{V}_n^{(5,1,1,1)}$  and $\mathcal{V}_n^{(1,0,4,5)}$, between $\mathcal{V}_n^{(3,2,8,3)}$ and $\mathcal{V}_n^{(3,1,4,3)}$.
\end{corollary}

\pf We first provide a bijection between $\mathcal{V}_n^{(4,3,7,2)}$ and $\mathcal{V}_n^{(2,1,7,4)}$, the others are similar and left to the interested readers. Note that $\mathcal{V}_n^{(4,3,7,2)}$ has the
weight sequences
$$\alpha_k=3^{k-1}, \beta_k=7\cdot 2^{k-1}, \gamma_k=7(3^{k-1}-2^{k-1})=7(2^{k-2}+2^{k-3}\cdot 3+\cdots+2\cdot 3^{k-3}+3^{k-2}), $$
and $\mathcal{V}_n^{(2,1,7,4)}$ has the weight sequences
$$\alpha'_k=1, \beta'_k=7\cdot 4^{k-1}, \gamma'_k=\frac{7}{3}(4^{k-1}-1)=7(4^{k-2}+4^{k-3}+\cdots+4+1), \ \ (k\geq 1). $$
Note that $\gamma_1=\gamma'_1=0$. One can give an equivalent weighting method for $\mathcal{V}_n^{(4,3,7,2)}$ and $\mathcal{V}_n^{(2,1,7,4)}$. For any primitive part $\mathbf{P}_1$ of any path $\mathbf{P}\in\mathcal{V}_n^{(4,3,7,2)}$, the weight of $\mathbf{P}_1$ is re-assigned as follows:

$1)$ If $\mathbf{P}_1=\mathbf{u}^k\mathbf{d}^k$, a maximal pyramid of height $k$ at altitude $0$, we first mark the endpoint of the $i$-th $\mathbf{u}$-step, then weight the $\mathbf{u}$-steps along the path respectively by $7, \underbrace{2, \dots , 2}_{i-1}, \underbrace{3, \dots , 3}_{j-1}, 1$ for $i+j=k$ with $i,j\geq 1$ and $k\geq 2$, the weight $2$ can also be regarded as $1$ or $\hat{1}$. Hence the total weight of $\mathbf{P}_1$ is just $\gamma_k$.

$2)$ If $\mathbf{P}_1=\mathbf{u}^k\mathbf{Q}\mathbf{d}^k$, where $\mathbf{Q}$ is a concatenation of at least two maximal pyramids at altitude $0$, we first mark the endpoint of the $k$-th $\mathbf{u}$-step of $\mathbf{P}_1$, then weight the $\mathbf{u}$-steps of any maximal pyramid of height $j$ along the path respectively by $\underbrace{3, \dots , 3}_{j-1}, 1$ for $j\geq 1$ and weight the first $k$ $\mathbf{u}$-steps along the path respectively by $7, \underbrace{2, \dots , 2}_{k-1}$ for $k\geq 1$, the weight $2$ can also be regarded as $1$ or $\hat{1}$.

For any primitive part $\mathbf{P}_1$ of any path $\mathbf{P}\in\mathcal{V}_n^{(2,1,7,4)}$, the weight of $\mathbf{P}_1$ is re-assigned as follows:

$1)$ If $\mathbf{P}_1=\mathbf{u}^k\mathbf{d}^k$, a maximal pyramid of height $k$ at altitude $0$, we first mark the endpoint of the $i$-th $\mathbf{u}$-step, then weight the $\mathbf{u}$-steps along the path respectively by $7, \underbrace{4, \dots , 4}_{i-1}, \underbrace{1, \dots , 1}_{j-1}, 1$ for $i+j=k$ with $i,j\geq 1$ and $k\geq 2$, the weight $4$ can also be regarded as $1$ or $\hat{3}$. Hence the total weight of $\mathbf{P}_1$ is just $\gamma'_k$.

$2)$ If $\mathbf{P}_1=\mathbf{u}^k\mathbf{Q}\mathbf{d}^k$, where $\mathbf{Q}$ is a concatenation of at least two maximal pyramids at altitude $0$, we first mark the endpoint of the $k$-th $\mathbf{u}$-step of $\mathbf{P}_1$, then weight each $\mathbf{u}$-step of any maximal pyramid by $1$ and weight the first $k$ $\mathbf{u}$-steps along the path respectively by $7, \underbrace{4, \cdots , 4}_{k-1}$ for $k\geq 1$, the weight $4$ can also be regarded as $1$ or $\hat{3}$.

Now we can give a recursive bijection $\tau$ between $\mathcal{V}_n^{(4,3,7,2)}$ and $\mathcal{V}_n^{(2,1,7,4)}$. For any $\mathbf{P}\in \mathcal{V}_n^{(4,3,7,2)}$, $\mathbf{P}$ can be uniquely partitioned into $\mathbf{P}=\mathbf{P}_1\mathbf{P}_2\dots \mathbf{P}_{\ell}$, where each $\mathbf{P}_h$ is a primitive weighted Dyck path, then define $\tau(\mathbf{P})=\tau(\mathbf{P}_1)\tau(\mathbf{P}_2)\dots \tau(\mathbf{P}_{\ell})\in \mathcal{V}_n^{(2,1,7,4)}$. So it is sufficed to consider that $\mathbf{P}$ is primitive.

Let $\mathbf{\mathbf{P}}=\mathbf{u}^{k_0}\mathbf{Q}\mathbf{d}^{k_0}\in \mathcal{V}_n^{(4,3,7,2)}$ be primitive and $\mathbf{Q}=(\mathbf{u}^{k_m}\mathbf{d}^{k_m})(\mathbf{u}^{k_{m-1}}\mathbf{d}^{k_{m-1}})\dots(\mathbf{u}^{k_1}\mathbf{d}^{k_1})$ with $m$ maximal pyramids such that $k_1+\cdots +k_m=n-k_0$ for $m, k_1, \dots, k_m\geq 1$ and $1\leq k_0<n$, $\mathbf{P}$ has a marked point at the end of the $k_0$-th $\mathbf{u}$-step and the $\mathbf{u}$-steps of $\mathbf{P}$ are weighted along the path by $7, w_1, \dots , w_{k_0-1}, v_1, v_2, \dots, v_{n-k_0}$, where $w_1, \dots, w_{k_0-1}\in \{1, \hat{1}\}$ and $(v_1, v_2, \dots, v_{n-k_0})=(~\underbrace{3, \dots , 3}_{k_m-1}, 1, \underbrace{3, \dots , 3}_{k_{m-1}-1}, 1, \dots, \underbrace{3, \dots , 3}_{k_1-1}, 1)$.
Write $w_1w_2\dots w_{k_0-1}=1^{s_1}\hat{1}^{r_1}\dots 1^{s_h}\hat{1}^{r_h}1^{s_{h+1}}$ such that $r_1+\cdots+r_h+s_1+\cdots +s_h+s_{h+1}=k_0-1$ for $h, s_1, s_{h+1}\geq 0$ and $r_1, \dots, r_h, s_2, \dots, s_{h}\geq 1$. Note that $\mathbf{Q}$ has another unique form, that is,
$$\mathbf{Q}=\big(\mathbf{u}^{i_p+1}\mathbf{d}^{i_p+1}(\mathbf{ud})^{j_p-1}\big)\big(\mathbf{u}^{i_{p-1}+1}\mathbf{d}^{i_{p-1}+1}(\mathbf{ud})^{j_{p-1}-1}\big)\dots \big(\mathbf{u}^{i_1+1}\mathbf{d}^{i_1+1}(\mathbf{ud})^{j_{1}-1}\big), $$
where $i_p\geq 0, p, i_1, \dots, i_{p-1}, j_1, \dots, j_p\geq 1$ and $i_1+\dots+i_p+j_1+ \dots +j_p=n-k_0$ with $i_p+1=k_m$. Then $\mathbf{Q}$ has the equivalent weight representation $v_1v_2\dots v_{n-k_0}=3^{i_p}1^{j_p}3^{i_{p-1}}1^{j_{p-1}}\dots 3^{i_1}1^{j_1}. $

\vskip0.2cm

Define
$$\tau(\mathbf{P})=\mathbf{u}^{1+(k_1-1)+k_2+\cdots+k_m}\mathbf{Q}'\mathbf{d}^{k_m+\cdots+k_2+(k_1-1)+1}, $$
where $\mathbf{Q}'=(\mathbf{ud})^{s_{h+1}}\big(\mathbf{u}^{r_h+1}\mathbf{d}^{r_h+1}(\mathbf{ud})^{s_{h}-1}\big)\dots \big(\mathbf{u}^{r_2+1}\mathbf{d}^{r_2+1}(\mathbf{ud})^{s_{2}-1}\big)\big(\mathbf{u}^{r_1+1}\mathbf{d}^{r_1+1}(\mathbf{ud})^{s_{1}}\big)$.
Replace the weight $3$ by $\hat{3}$ and $\hat{1}$ by $1$, and mark $\tau(\mathbf{P})$ at the endpoint of the $(n-k_0)$-th $\mathbf{u}$-step. Note that the $\mathbf{u}$-steps of $\tau(\mathbf{P})$ are weighted along the path by $7, \overline{v}_1, \overline{v}_2, \dots, \overline{v}_{n-k_0-1}, \underbrace{1, 1, \dots, 1}_{k_0}$, where $\overline{v}_1\overline{v}_2\dots\overline{v}_{n-k_0-1}=1^{j_1-1}\hat{3}^{i_1}1^{j_2}\hat{3}^{i_2}\dots 1^{j_p}\hat{3}^{i_p}$. Hence, $\tau(\mathbf{P})\in \mathcal{V}_n^{(2,1,7,4)}$.

Conversely, let $\mathbf{\mathbf{P}'}=\mathbf{u}^{n-k_0}\mathbf{Q}''\mathbf{d}^{n-k_0}\in \mathcal{V}_n^{(2,1,7,4)}$ be primitive and $$\mathbf{Q}''=(\mathbf{u}^{q_m}\mathbf{d}^{q_m})(\mathbf{u}^{q_{m-1}}\mathbf{d}^{q_{m-1}})\dots(\mathbf{u}^{q_1}\mathbf{d}^{q_1})$$
with $m$ maximal pyramids such that $q_1+\cdots +q_m=k_0$ for $m, q_1, \dots, q_m\geq 1$ and $1\leq k_0<n$, $\mathbf{P}'$ has a marked point at the end of the $(n-k_0)$-th $\mathbf{u}$-step and the $\mathbf{u}$-steps of $\mathbf{P}'$ are weighted along the path by $7, \overline{\overline{v}}_1, \overline{\overline{v}}_2, \dots, \overline{\overline{v}}_{n-k_0-1}, \overline{w}_1, \overline{w}_2, \dots, \overline{w}_{k_0}$, where $\overline{\overline{v}}_1\overline{\overline{v}}_2\dots\overline{v}_{n-k_0-1}=1^{j_1-1}\hat{3}^{i_1}1^{j_2}\hat{3}^{i_2}\dots 1^{j_p}\hat{3}^{i_p}$ for $i_p\geq 0, p, i_1, \dots, i_{p-1}, j_1, \dots, j_p\geq 1$ and $i_1+\dots+i_p+j_1+ \dots +j_p=n-k_0$, and $(\overline{w}_1, \overline{w}_2, \dots, \overline{w}_{k_0})=(\underbrace{1, \dots , 1}_{q_m-1}, 1, \underbrace{1, \dots , 1}_{q_{m-1}-1}, 1, \dots, \underbrace{1, \dots , 1}_{q_1-1}, 1)$ or equivalently $(\overline{w}_1, \overline{w}_2, \dots, \overline{w}_{k_0})=(\underbrace{\hat{1}, \dots , \hat{1}}_{q_m-1}, 1, \underbrace{\hat{1}, \dots , \hat{1}}_{q_{m-1}-1}, 1, \dots, \underbrace{\hat{1}, \dots , \hat{1}}_{q_1-1}, 1)$. Write
$$\overline{w}_1\overline{w}_2\dots \overline{w}_{k_0}=1^{s_{h+1}}\hat{1}^{r_h}1^{s_h}\dots \hat{1}^{r_1}1^{s_1}$$
such that $r_1+\cdots+r_h+s_1+\cdots +s_h+s_{h+1}=k_0$ for $h, s_{h+1}\geq 0$ and $r_1, \dots, r_h, s_1, \dots, s_{h}\geq 1$. Define
$$\tau'(\mathbf{P}')=\mathbf{u}^{1+(q_1-1)+q_2+\cdots+q_m}\mathbf{Q}'''\mathbf{d}^{q_m+\cdots+q_2+(q_1-1)+1}, $$
where $\mathbf{Q}'''=\big(\mathbf{u}^{i_p+1}\mathbf{d}^{i_p+1}(\mathbf{ud})^{j_{p}-1}\big)\dots \big(\mathbf{u}^{i_2+1}\mathbf{d}^{i_2+1}(\mathbf{ud})^{j_{2}-1}\big)\big(\mathbf{u}^{i_1+1}\mathbf{d}^{i_1+1}(\mathbf{ud})^{j_{1}-1}\big)$.
Replace the weight $\hat{3}$ by $3$, and mark $\tau'(\mathbf{P}')$ at the endpoint of the $k_0$-th $\mathbf{u}$-step. Note that the $\mathbf{u}$-steps of $\tau'(\mathbf{P}')$ are weighted along the path by $7,  \overline{\overline{w}}_1, \overline{\overline{w}}_2, \dots, \overline{\overline{w}}_{k_0-1}, v_1', v_2', \dots, v_{n-k_0}'$, where $\overline{\overline{w}}_1\overline{\overline{w}}_2\dots\overline{\overline{w}}_{k_0-1}= 1^{s_1-1}\hat{1}^{r_1}\dots 1^{s_h}\hat{1}^{r_h}1^{s_{h+1}}$, $v_1' v_2'\dots v_{n-k_0}'={3}^{i_p}1^{j_p}\dots {3}^{i_2}1^{j_2}{3}^{i_1}1^{j_1}$. Hence, this shows that $\tau'(\mathbf{P}')\in \mathcal{V}_n^{(4,3,7,2)}$.

One can easily testify that $\tau(\tau'(\mathbf{P}'))=\mathbf{P}'$ and $\tau'(\tau(\mathbf{P}))=\mathbf{P}$, namely, $\tau$ is indeed a bijection between $\mathcal{V}_n^{(4,3,7,2)}$ and $\mathcal{V}_n^{(2,1,7,4)}$. \qed

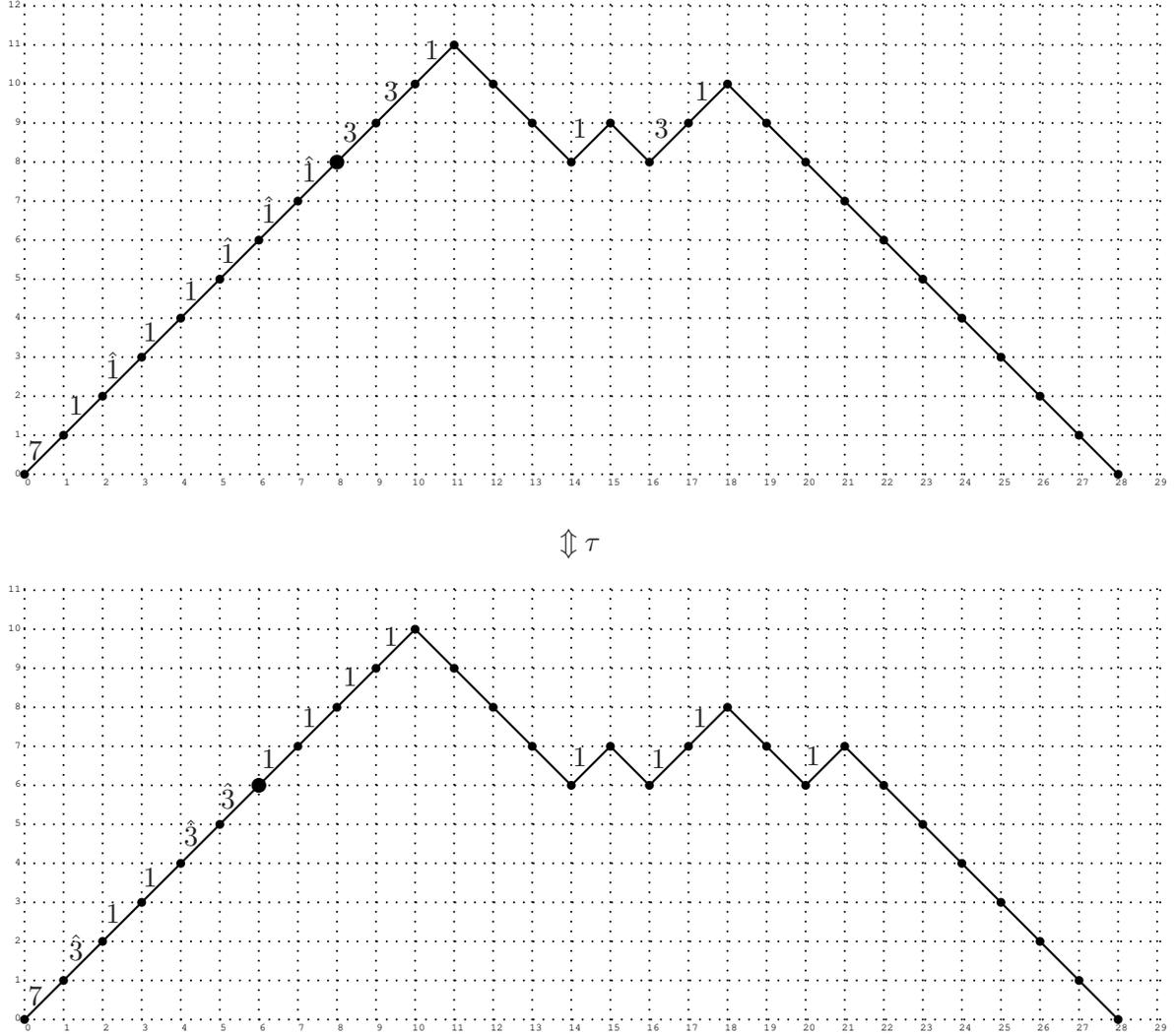
\begin{figure}[h] \setlength{\unitlength}{0.5mm}

\begin{center}
\begin{pspicture}(15,6.5)
\psset{xunit=15pt,yunit=15pt}\psgrid[subgriddiv=1,griddots=4,
gridlabels=4pt](0,0)(29,12)

\psline(0,0)(8,8)(11,11)(14,8)(15, 9)(16,8)(18, 10)(28,0)

\pscircle*(0,0){0.06}\pscircle*(1,1){0.06}\pscircle*(2,2){0.06}
\pscircle*(3,3){0.06}\pscircle*(4,4){0.06}\pscircle*(5,5){0.06}
\pscircle*(6,6){0.06}\pscircle*(7,7){0.06}\pscircle*(8,8){0.1}
\pscircle*(9,9){0.06}\pscircle*(10,10){0.06}\pscircle*(11,11){0.06}
\pscircle*(12,10){0.06}\pscircle*(13,9){0.06}\pscircle*(14,8){0.06}
\pscircle*(15,9){0.06}\pscircle*(16,8){0.06}\pscircle*(17,9){0.06}
\pscircle*(18,10){0.06}\pscircle*(19,9){0.06}

\pscircle*(28,0){0.06}\pscircle*(27,1){0.06}\pscircle*(26,2){0.06}
\pscircle*(25,3){0.06}\pscircle*(24,4){0.06}\pscircle*(23,5){0.06}
\pscircle*(22,6){0.06}\pscircle*(21,7){0.06}\pscircle*(20,8){0.06}

\put(0.05,.2){$7$}\put(.6,.8){$1$}\put(1.1,1.3){$\hat{1}$}\put(1.6,1.8){$1$}\put(2.15,2.35){$1$}\put(2.65,2.85){$\hat{1}$}
\put(3.2,3.4){$\hat{1}$}\put(3.75,3.95){$\hat{1}$}

\put(4.3,4.5){$3$}\put(4.85,5.05){$3$}\put(5.4,5.6){$1$} \put(7.4,4.55){$1$}\put(8.5,4.55){$3$}\put(9.05,5.05){$1$}

\end{pspicture}
\end{center}\vskip0.5cm

$\Updownarrow \tau$
\vskip0.7cm

\begin{center}
\begin{pspicture}(15,5.5)
\psset{xunit=15pt,yunit=15pt}\psgrid[subgriddiv=1,griddots=4,
gridlabels=4pt](0,0)(29,11)

\psline(0,0)(6,6)(10,10)(14,6)(15, 7)(16,6)(18, 8)(20,6)(21,7)(22,6)(28,0)

\pscircle*(0,0){0.06}\pscircle*(1,1){0.06}\pscircle*(2,2){0.06}
\pscircle*(3,3){0.06}\pscircle*(4,4){0.06}\pscircle*(5,5){0.06}
\pscircle*(6,6){0.1}\pscircle*(7,7){0.06}\pscircle*(8,8){0.06}
\pscircle*(9,9){0.06}\pscircle*(10,10){0.06}\pscircle*(11,9){0.06}
\pscircle*(12,8){0.06}\pscircle*(13,7){0.06}\pscircle*(14,6){0.06}
\pscircle*(15,7){0.06}\pscircle*(16,6){0.06}\pscircle*(17,7){0.06}
\pscircle*(18,8){0.06}\pscircle*(19,7){0.06}\pscircle*(20,6){0.06}

\pscircle*(28,0){0.06}\pscircle*(27,1){0.06}\pscircle*(26,2){0.06}
\pscircle*(25,3){0.06}\pscircle*(24,4){0.06}\pscircle*(23,5){0.06}
\pscircle*(22,6){0.06}\pscircle*(21,7){0.06}

\put(0.05,.2){$7$}\put(.6,.8){$\hat{3}$}\put(1.1,1.3){$1$}\put(1.6,1.8){$1$}\put(2.15,2.35){$\hat{3}$}\put(2.65,2.85){$\hat{3}$}
\put(3.2,3.4){$1$}\put(3.75,3.95){$1$}\put(4.3,4.5){$1$}\put(4.85,5.05){$1$}

\put(7.4,3.42){$1$}\put(8.46,3.4){$1$}\put(9.02,3.95){$1$}\put(10.55,3.43){$1$}

\end{pspicture}
\end{center}

\caption{\small An example of the bijection $\tau$ described in the proof of Corollary 3.12. }

\end{figure}

In order to give a more intuitive view on the bijection $\tau$, we present a pictorial description of
$\tau$ for the case $\mathbf{P}=\mathbf{u}^8(\mathbf{u}^3\mathbf{d}^3)(\mathbf{u}\mathbf{d})(\mathbf{u}^2\mathbf{d}^2)\mathbf{d}^8\in \mathcal{V}_n^{(4,3,7,2)}$, $\mathbf{P}$ has a marked point at the end of the $8$-th $\mathbf{u}$-step and its $\mathbf{u}$-steps are weighted along the path by $7, 1, \hat{1}, 1, 1, \hat{1}, \hat{1},\hat{1}, 3, 3, 1, 1, 3, 1$, and $\tau(\mathbf{P})=\mathbf{u}^6(\mathbf{u}^4\mathbf{d}^4)(\mathbf{u}\mathbf{d})(\mathbf{u}^2\mathbf{d}^2)(\mathbf{u}\mathbf{d})\mathbf{d}^6\in \mathcal{V}_n^{(2,1,7,4)}$, $\tau(\mathbf{P})$ has a marked point at the end of the $6$-th $\mathbf{u}$-step and its $\mathbf{u}$-steps are weighted along the path by $7, \hat{3}, 1, 1, \hat{3}, \hat{3},1,1,1,1,1,1,1, 1$. See Figure 5.

\begin{remark}
Note that for $n\geq 2$, there hold
\begin{eqnarray*}
\sum_{i=0}^{n-2}D_iD_{n-2-i} &=& \frac{1}{7}|\mathcal{V}_n^{(4,3,7,2)}|=\frac{1}{7}|\mathcal{V}_n^{(2,1,7,4)}| = \frac{1}{4}|\mathcal{V}_n^{(5,4,4,1)}|   \\
                             &=& \frac{1}{4}|\mathcal{V}_n^{(5,1,1,1)}|=\frac{1}{4}|\mathcal{V}_n^{(1,0,4,5)}|=\frac{1}{8}|\mathcal{V}_n^{(3,2,8,3)}|=\frac{1}{8}|\mathcal{V}_n^{(3,1,4,3)}|.
\end{eqnarray*}
One can be asked whether there exist direct combinatorial interpretations by Delannoy paths.
\end{remark}

\subsection{The special case related to complete $(r+1)$-ary trees}

Let $T(x)=\sum_{n\geq 0}T_nx^n$ be the generating function for Fuss-Catalan numbers $T_n=\frac{1}{nr+1}\binom{n(r+1)}{n}$, which count the number of complete $(r+1)$-ary trees with $n$ internal vertices \cite{StanleyEC}. $T(x)$ obeys the relations $T(x)=1+xT(x)^{r+1}=\frac{1}{1-xT(x)^{r}}=\frac{1}{1-xT(x)^{r-1}-x^2T(x)^{2r}}$.

\begin{example} The special case when
$$\alpha(x)=\beta(x)=xT(x)^{m}, \gamma(x)=x^2T(x)^{2m} $$
produces
\begin{eqnarray*}
V_{\alpha, \beta}(x)=1+\frac{\alpha(x)^2}{1-\alpha(x)-\alpha(x)^2}=1+\sum_{k=1}^{\infty}F_{k}\cdot\big(xT(x)^{m}\big)^{k+1}.
\end{eqnarray*}
By Lagrange's inversion formula \cite{Comtet, Gessel}, we have
$$V_n=\sum_{k=1}^{n-1}\frac{m(k+1)F_k}{n(r+1)+(m-r-1)(k+1)}\binom{n(r+1)+(m-r-1)(k+1)}{n-k-1},\ (n\geq 2),$$
with $V_0=1, V_1=0.$

The special case when
$$\alpha(x)=xT(x)^{r}, \beta(x)=xT(x)^{m}, \gamma(x)=x^2T(x)^{m+r} $$
produces
\begin{eqnarray*}
V_{\alpha, \beta}(x)=\frac{1-xT(x)^{r}}{1-xT(x)^{r}-x^2T(x)^{m+r}}=\frac{1}{1-(T(x)-1)^{2}T(x)^{m-r-1}}.
\end{eqnarray*}
By Lagrange's inversion formula \cite{Comtet, Gessel}, we have
$$V_n=\sum_{k=0}^{[\frac{n}{2}]}\sum_{j=0}^{n-2k}\binom{k(m-r-1)}{j}\frac{2k+j}{n}\binom{n(r+1)}{n-2k-j},\ (n\geq 1),$$
with $V_0=1.$ The $m=r+1$ case collapses to
$$V_n=\sum_{k=0}^{[\frac{n}{2}]}\frac{2k}{n}\binom{n(r+1)}{n-2k},\ (n\geq 1).$$
\end{example}

\begin{remark}
One can also consider some examples that $\gamma(x)\neq \alpha(x)\beta(x)$. For example, the special case when
$$\alpha(x)=xT(x)^{r}, \beta(x)=xT(x)^{m}, \gamma(x)=xT(x)^{r} $$
produces
\begin{eqnarray*}
V_{\alpha, \beta, \gamma}(x)=\frac{1}{1-xT(x)^{r}-\frac{x^3T(x)^{m+2r}}{1-xT(x)^{r}}}=\frac{T(x)}{1-(T(x)-1)^{3}T(x)^{m-r-1}}.
\end{eqnarray*}
By Lagrange's inversion formula, we have
$$V_n=\sum_{k=0}^{[\frac{n}{3}]}\sum_{j=0}^{n-3k}\binom{k(m-r-1)+1}{j}\frac{3k+j}{n}\binom{n(r+1)}{n-3k-j},\ (n\geq 1),$$
with $V_0=1.$ The $m=r+1$ case collapses to
$$V_n=\sum_{k=0}^{[\frac{n}{3}]} \frac{3kr+3k+1}{nr+3k+1}\binom{n(r+1)}{n-3k},\ (n\geq 1).$$

\end{remark}

\vskip.2cm

\section*{Declaration of competing interest}

The authors declare that they have no known competing financial interests or personal relationships that could have
appeared to influence the work reported in this paper.

\section*{Acknowledgements} {The authors are grateful to the referees for
the helpful suggestions and comments. The Project is sponsored by ``Liaoning
BaiQianWan Talents Program" and by the Fundamental Research Funds for the Central Universities Under Contract. }


\end{document}